\documentstyle[11pt]{article}
\title{\LARGE Improving $L^2$ estimates to Harnack inequalities }

\author{\Large Stathis Filippas$^{1,4}$ \& Luisa Moschini$^{2}$
\& Achilles Tertikas$^{3,4}$  \\
                                                                           \\
        Department of Applied Mathematics$^{1}$ \\
         University of Crete,
         71409 Heraklion,  Greece \\
        filippas@tem.uoc.gr\\
                                          \\
                                          Dipartimento di Metodi e Modelli Matematici
                                          per le Scienze Applicate$^2$,\\
 University of Rome ``La Sapienza'', 00161 Rome, Italy \\
moschini@dmmm.uniroma1.it \\
\\
 Department of Mathematics$^{3}$ \\
         University of Crete,
         71409 Heraklion,  Greece \\
          tertikas@math.uoc.gr\\                                      \\
        Institute of Applied and Computational Mathematics$^4$, \\
        FORTH, 71110 Heraklion, Greece \\
    \\ }

\begin{document}

\maketitle

\def\finedim{$\hfill \Box$}


\newcommand{\be}{\begin{equation}}
\newcommand{\ee}{\end{equation}}
\newcommand{\bea}{\begin{eqnarray}}
\newcommand{\eea}{\end{eqnarray}}
\newcommand{\la}{\label}
\newcommand{\xa}{\alpha}
\newcommand{\xb}{\beta}
\newcommand{\xg}{\gamma}
\newcommand{\xG}{\Gamma}
\newcommand{\xd}{\delta}
\newcommand{\xD}{\Delta}
\newcommand{\xe}{\varepsilon}
\newcommand{\xz}{\zeta}
\newcommand{\xh}{\eta}
\newcommand{\th}{\theta}
\newcommand{\Th}{\Theta}
\newcommand{\xk}{\kappa}
\newcommand{\xl}{\lambda}
\newcommand{\xL}{\Lambda}
\newcommand{\xm}{\mu}
\newcommand{\xn}{\nu}
\newcommand{\ks}{\xi}
\newcommand{\KS}{\Xi}
\newcommand{\xp}{\pi}
\newcommand{\xP}{\Pi}
\newcommand{\xr}{\rho}
\newcommand{\xs}{\sigma}
\newcommand{\xS}{\Sigma}
\newcommand{\xf}{\phi}
\newcommand{\xF}{\Phi}
\newcommand{\ps}{\psi}
\newcommand{\PS}{\Psi}
\newcommand{\xo}{\omega}
\newcommand{\xO}{\Omega}
\newcommand{\Ren}{ I \! \! R^N}
\newcommand{\R}{ I \! \! R}
\newcommand{\Ded}{{\cal D}^{1,2}}
\newcommand{\ho}{H_{\xO}}
\newcommand{\co}{C_{\xO}}
\newcommand{\k}{{\biggl(\frac{N-2}{2}\biggr)}^2}
\newcommand{\weight}{|x|^{-(N-2)}}
\newcommand{\Lp}{L^p (\xO)}
\newcommand{\spc}{{\cal C}_0^{\infty} (\xO)}
\newcommand{\sph}{H_0^1 (\xO)}
\newcommand{\spw}{W_0^{1,2} (\xO ;\;  \weight )}
\newcommand{\iuu}{\int_{\Ren} |\nabla u|^2\,dx}
\newcommand{\ioo}{\int_{\Ren} Vu^2\,dx}
\newcommand{\Br}{B_r}
\newcommand{\bBr}{\partial \! B_r}
\newcommand{\ra}{\rightarrow}
\newcommand{\rft}{\rightarrow +\infty}
\newcommand{\bin}{\int_{\bBr}}
\newcounter{newsection}
\newtheorem{theorem}{Theorem}[section]
\newtheorem{lemma}[theorem]{Lemma}
\newtheorem{proposition}[theorem]{Proposition}
\newtheorem{corollary}[theorem]{Corollary}
\newtheorem{remark}[theorem]{Remark}
\newtheorem{definition}[theorem]{Definition}
\newcounter{newsec} \renewcommand{\theequation}{\thesection.\arabic{equation}}

\begin{abstract}
We consider operators of the form ${\mathcal L}=-L-V$, where $L$
is an elliptic operator and $V$ is a singular potential, defined
on a smooth bounded domain $\Omega\subset \R^n$ with Dirichlet
boundary conditions. We allow the boundary of $\Omega$ to be made
of various pieces of different codimension. We assume that
${\mathcal L}$ has a generalized first eigenfunction of which we
know two sided estimates. Under these assumptions we prove
optimal Sobolev inequalities for the operator ${\mathcal L}$, we
show that it generates an intrinsic ultracontractive semigroup and
finally we derive a parabolic Harnack inequality up to the
boundary as well as sharp heat kernel estimates.
\smallskip

\noindent {\bf AMS Subject Classification: }35K65, 26D10 (35K20, 35B05.)  \\
{\bf Keywords: } Singular heat equation, optimal Sobolev
inequalities, generalized first eigenvalue, parabolic Harnack
inequality, heat kernel estimates, degenerate elliptic operators.
\end{abstract}

\setcounter{equation}{0}
\section{Introduction}

Let $\Omega \subset \R^n$,  be a bounded domain and suppose $V$ is
an $L^1_{loc}(\Omega)$ potential for which we know the following
$L^2$ estimate \be \label{ref1} 0 < \xl_1 := \inf_{u \in
C_0^{\infty}(\xO)} \frac{ \int_{\Omega} (|\nabla u|^2-V(x) u^2) dx }
{ \int_{\xO} u^{2} dx } \ .\ee One of the motivations of the present
work is whether one can improve the above estimate to a Sobolev type
estimate, involving, if possible, the critical Sobolev exponent. It
is clear that to improve the previous estimate one needs more
information concerning the potential $V$, besides (\ref{ref1}). One
additional piece of information that we are going to use,  is the
existence of a generalized eigenfunction $\phi_1$ of problem
(\ref{ref1}) as well as sharp  two sided estimates
 of $\phi_1$. Under this extra piece of information we are
able to obtain sharp Sobolev type inequalities involving the
critical Sobolev exponent.

The knowledge of the asymptotic behaviour of $\phi_1$ usually comes
as a consequence of the maximum principle and the local character of
$V$. We present such an argument, following the ideas of Brezis,
Marcus and Shafrir  \cite{BMS} for the critical potential
$V(x)=\frac{1}{4} \frac{1}{{\rm dist}^2 (x,\partial \Omega)}$ in the
Appendix.   We should mention that in this case the asymptotics of
$\phi_1$ have already been derived  by D\'{a}vila and Dupaigne
\cite{DD1}, \cite{DD2}. The argument  we present  is much simpler
and is based only on the maximum principle applied in the
appropriate energy space.
 All the potentials we have in Section 4
as well as other potentials can be treated similarly.

Presupposing the existence and asymptotic behaviour of the
generalized eigenfunction $\phi_1$ seems to be a natural
assumption. It is E. B Davies  \cite{D2}  who put forward the idea
of  connecting the asymptotics of  $\phi_1$ to the asymptotics
of the Green function  in the
case of subcritical potentials  $V$.
 In fact, he conjectured that
knowing that the asymptotic behavior of $\phi_1$ is like
$d(x):={\rm dist} (x,\partial \Omega)$ is actually equivalent to
the exact two-sided Green function bounds. This conjecture was
answered positively in \cite{FMT1} even for a larger class of
potentials for which the generalized eigenfunction $\phi_1$
behaves like $d^\alpha(x):={\rm dist}^\alpha (x,\partial \Omega)$,
for some $\alpha \ge 1/2$.

The idea of obtaining heat kernel estimates of second order
elliptic operators with singular potentials in terms of the
generalized ground state is not new and besides  \cite{G} and
\cite{GSC} has been successfully exploited in  \cite{LS} and
\cite{MS}.

Another motivation is the study of the corresponding parabolic
problem, especially when the potential $V$ is singular (see e.g.
\cite{BG}, \cite{CM} and \cite{VZ}). In connection with this, we
mention the work of Cabr\'e and Martel \cite{CM} where the
condition $\lambda_1>-\infty$ is shown to be necessary and
sufficient for the existence of a global positive weak solution.
They also show that these solutions grow at most exponentially in
time for any nonnegative initial data $u_0\in L^2(\Omega)$. In
fact, as far as the parabolic problem is concerned, condition
(\ref{ref1}) is practically equivalent to $\lambda_1>-\infty$ due
to a shift in the time variable. Under the extra assumption that
the asymptotics of the generalized eigenfunction are known, we
show that the corresponding Schr\"odinger operator generates a
semigroup of integral operators for which we obtain parabolic
Harnack inequality up to the boundary and precise heat kernel
estimates.

At this point we introduce some notations that we keep throughout
the work. We assume that $\xO \subset \R^n$, $n \geq 2$, is a
bounded domain, with a boundary $\partial \xO = \cup_{k=1}^{n}
\xG_{k}$, where
 $ \xG_{k} = \cup_{j=1}^{m_k}  \xG_{k,j}  $  is  a finite
union    of $m_k$   distinct  smooth $C^2$  boundaryless
hypersurfaces $\xG_{k,j} $
 of codimension ${\rm codim} \xG_{k,j}
=k$,  where   $j=1,\ldots,m_k$,        $k=1,\ldots,n-1$;
  in addition $\xG_n = \{x_1,
\dots, x_{m_n} \}$,   $\xG_{k,j} \cap \xG_{l,i} = \emptyset$  if $k
\neq l$, or $i \neq j$  and $\xG_{k,j} \cap \xG_n= \emptyset$, for
$k=1, \ldots, n-1$ and $j=1,\ldots,m_k$. We also set $d_k(x) =  {\rm
dist}(x, \xG_k)$. For $x\in \Omega$  we  denote by $d(x)$ the
distance to the boundary $\partial \xO$. We clearly have that $d(x)=
{\rm min}_{x \in \xO} \{ d_1(x), \ldots, d_n(x) \}$. Finally for the
part of the boundary that is of codimension one we use the special
notation $\partial_1 \Omega=\Gamma_1$.

We are interested in the quadratic form \be\la{3.1} Q[u] =
\int_{\xO} \left( \sum_{i,j=1}^{n} a_{ij}(x) \frac{\partial
u}{\partial x_i}\frac{\partial u}{\partial x_j} - V u^2 \right)
dx,~~~~~~~~u \in C_0^{\infty}(\xO), \ee where $V \in
L^1_{loc}(\xO)$ and $ a_{ij}(x)$, $i,j=1,\ldots,n$ is a measurable
symmetric uniformly elliptic matrix, that is \be\la{3.ell} C_0
|\xi|^2 \leq \sum_{i,j=1}^{n} a_{ij}(x) \xi_i \xi_j  \leq C_0^{-1}
|\xi|^2, \ \ ~~~~~~ \xi \in \R^n \ee for some $C_0>0$. We also
assume the $L^2$ estimate \be\la{3.5} 0 < \xl_1 := \inf_{u \in
C_0^{\infty}(\xO)} \frac{ Q[u] } { \int_{\xO} u^{2} dx }, \ee and
that to $\xl_1$, there
 corresponds a generalized eigenfunction $\phi_1$
 of (\ref{3.5}). More precisely, we assume
that $\phi_1  \in H^1_{loc}(\xO)  \cap  L^{\infty}_{loc}(\xO)$ and
that \be\la{1.10b}
 \int_{\xO}  \sum_{i,j=1}^{n} a_{ij} \frac{\partial \phi_1}{\partial x_i}  \frac{\partial \psi}{\partial x_j}
 \,  dx  =
\int_{\xO} (V + \xl_1) \phi_1 \psi \,  dx, ~~~~~~~~  \ \psi
\in C_0^{\infty}(\xO). \ee
 In addition we assume that we have  two sided estimates on
  $\phi_1$  of the form
 $\phi_1  \sim  d_1^{\xa_1} \ldots d_n^{\xa_n}  $, that is \be\la{2.3x}  c_1
d_1^{\xa_1}(x) \ldots d_n^{\xa_n}(x) \leq \phi_1(x)  \leq c_2
d_1^{\xa_1}(x) \ldots d_n^{\xa_n}(x), \ee for any $x\in \Omega$, for
two positive constants $c_1,c_2$ and for suitable exponents
$\alpha_1,\ldots,\alpha_n$. Appropriate conditions on the exponents
$\alpha_i$ will be formulated below.

Our first result concerns the following improved Sobolev type
inequality.

\begin{theorem}{\bf (Optimal Sobolev type inequality)} \la{thm3.1zero}
For  $V \in  L^1_{loc}(\xO)$ we assume that (\ref{3.5}) holds and
in addition there exists a  ground state $\phi_1  \in
H^1_{loc}(\xO)  \cap  L^{\infty}_{loc}(\xO)$ satisfying
(\ref{2.3x}) for \be\la{3.8zero} \xa_k > -\frac{k-2}{2} - (n-k)
\frac{q-2}{2(q+2)} \ ,  \ \ k=1,\ldots,n \ ,\ee where $2 <  q \leq
\frac{2n}{n-2}$ if $n\ge 3$ and $q>2$ if $ n=2 .$
 We then  have that
\be\la{3.9zero} 0 < C(\xO,\alpha_1,\ldots, \alpha_n, q) =
 \inf_{u \in C_0^{\infty}(\xO)} \frac{Q[u]}
{ \left( \int_{\xO}d^{\frac{q(n-2)-2n}{2}}
 |u|^{q} dx \right)^{\frac{2}{q}}} \ .
\ee In particular when $n \ge 3$ and $q= \frac{2n}{n-2}$
we have that  
$$ 0 < C(\xO,\alpha_1,\ldots,\alpha_n) =
 \inf_{u \in C_0^{\infty}(\xO)} \frac{Q[u]}
{ \left( \int_{\xO}
 |u|^{\frac{2n}{n-2}} dx \right)^{\frac{n-2}{n}}}
$$ if 
$$\xa_k > -\frac{k-2}{2} -
\frac{(n-k)}{2(n-1)},~~~~~~~~k=1,\ldots,n \ . $$
\end{theorem}

We note that the condition (\ref{3.8zero}) is optimal. In the
extreme cases $q>2$ and $\alpha_n=-\frac{n-2}{2}$  or $q=2$ and
$\alpha_k=-\frac{k-2}{2}$ we have different improved inequalities,
see Theorems \ref{thm3.2} and \ref{th2.2}  respectively.

Using the above Sobolev type inequality we will now proceed to the
study of the corresponding parabolic problem, that is \be
\label{defl} \left\{\begin{array}{lll} \frac{\partial u}{\partial
t}=-{\mathcal L} u:=\sum^{n}_{i,j=1} \frac{\partial}{\partial x_i}
\left(a_{ij}(x) \frac{\partial
u}{\partial x_j}\right)+ V(x) u  & \hbox { in } \ \ (0,\infty)\times \Omega ~ ,\\
u(x,t)=0 & \ \hbox
{ on } \ \ (0,\infty)\times \partial_1\Omega , \\
u(x,0)=u_0(x)  & \ \hbox { on } \ \ \Omega \ .\end{array} \right.
\ \ \ee

Our first result is the following Harnack inequality.

\begin{theorem}{\bf (Parabolic Harnack inequality up to the
boundary)} \label{harnackgen2} For  $V \in  L^1_{loc}(\xO)$ we
assume (\ref{3.5}) and that (\ref{2.3x}) holds for some $\alpha_k
\ge -\frac{k-2}{2}$, for $k=1,\cdots, n$. Then for $\mathcal L$ as
in (\ref{defl}) there exist positive constants $C_H$ and
$R=R(\Omega)$ such that for $x\in \Omega$, $0<r<R$ and for any
positive solution $u(y,t)$ of
$$\frac{\partial u}{\partial
t}=-{\mathcal L } u \hbox{ in } \left\{{\mathcal B}(x,r)\cap
\Omega\right\}\times (0,r^2) \ ,
$$ the following estimate holds true
$${\rm ess~sup}_{(y,t)\in \left\{{\mathcal B}(x,\frac{r}{2})\cap \Omega
\right\}\times (\frac{r^2}{4},\frac{r^2}{2})} \
u(y,t)\prod^n_{i=1} d_i^{-\alpha_i}(y) \le$$ $$\le C_H ~ {\rm
ess~inf}_{(y,t)\in \left\{{\mathcal B}(x,\frac{r}{2})\cap \Omega
\right\}\times (\frac{3}{4} r^2,r^2)}\  u(y,t) \prod^n_{i=1}
d_i^{-\alpha_i}(y) .$$
\end{theorem}

 Here ${\mathcal B}(x,r)$ denotes roughly speaking an $n$
dimensional cube centered at $x$ and having size $r$, see Definition
\ref{ball}  for details.

Theorem \ref{harnackgen2} states a parabolic Harnack inequality up
to the boundary for the ratio of any positive local solution to the
Cauchy-Dirichlet problem and the generalized eigenfunction $\phi_1$.
We note that $\alpha_1\ge 1/2$, therefore $\phi_1$ is zero on the
boundary $\partial_1 \Omega$. In particular it implies that any two
nonnegative solutions vanishing on $\partial_1 \Omega$ must vanish
at the same rate. It is clear then that such a normalization is
necessary. In fact the natural quantity is $v=u / \phi_1$ and it is
for this function that we prove the Harnack inequality. For the definition
of solution for   the function $v$ we refer to  Definition \ref{sol},
where however the appropriate weight is $\phi_1^2$.

 Alternatively, one
could define  local  weak solutions of (\ref{defl}) using  a suitable local  energetic
space obtained via  the quadratic form (\ref{3.1}). For an example of  globally defined  energetic
solutions see \cite{VZ}.

 Our result in the case $\alpha_1=1$, $\alpha_k=0$ for $k=2,
\cdots, n$ is basically the local comparison principle of
\cite{FGS}.
We also note that the restriction on $\alpha_k$ in
Theorem \ref{harnackgen2} is sharp.

In what follows  we denote by $h$ the integral kernel of the $L^2$
semigroup associated to the elliptic operator ${\mathcal L}$ as
defined in (\ref{defl}), that is
$$u(x,t):=\int_{\Omega} h(t,x,y) u_0(y) dy \ .$$ The existence of $h(t,x,y)$ is proved in
Proposition \ref{intrins} and it is a consequence of our Theorem
\ref{thm3.1zero}.

As usual, from the parabolic Harnack inequality one can
obtain sharp heat kernel estimates,  as explained by Grigoryan and Saloff Coste,
  see  \cite{G}, \cite{GSC}, \cite{SC}. In particular we have

\begin{theorem}
\label{heatgen2}{\bf (Sharp heat kernel estimates)}  For  $V \in
L^1_{loc}(\xO)$ we assume (\ref{3.5}) and that (\ref{2.3x}) holds
for some $\alpha_k \ge -\frac{k-2}{2}$, for $k=1,\cdots, n$.  Then
there exist positive constants $C_1, C_2$, with $C_1\le C_2$, and
$T>0$ depending on $\Omega$ such that
$$C_1 \prod^{n}_{i=1} \left(1+\frac{\sqrt t}{d_i(x)}\right)^{-\alpha_i}
\left(1+\frac{\sqrt t}{d_i(y)}\right)^{-\alpha_i} t^{-\frac{n}{2}}
e^{-C_2\frac{|x-y|^2}{t}}\le h(t,x,y) \le
$$ $$\le C_2 \prod^{n}_{i=1} \left(1+\frac{\sqrt t}{d_i(x)}\right)^{-\alpha_i}
\left(1+\frac{\sqrt t}{d_i(y)}\right)^{-\alpha_i} t^{-\frac{n}{2}}
e^{-C_1\frac{|x-y|^2}{t}}$$ for all $x,y\in \Omega$ and $0<t\le
T$, whereas
$$C_1 \prod^{n}_{i=1} d_i^{\alpha_i}(x)\
d_i^{\alpha_i}(y)\  e^{-\lambda_1 t} \le h(t,x,y) \le C_2
\prod^{n}_{i=1}d_i^{\alpha_i}(x)\ d_i^{\alpha_i}(y) \
e^{-\lambda_1 t}$$ for all $x,y\in \Omega$ and $t\ge T$.
\end{theorem}

Due to a shift in the time variable, Theorems \ref{harnackgen2} and
\ref{heatgen2} remain valid if we replace assumption (\ref{3.5})
with the condition $\lambda_1>-\infty$. For the corresponding
statement of Theorem \ref{thm3.1zero} under the condition
$\lambda_1>-\infty$ we refer to Theorem \ref{thm3.1}.

Although we present here only heat kernel estimates one can
integrate in time to obtain the corresponding Green function
estimates provided that $\lambda_1>0$.

It is clear that the asymptotics of $\phi_1$ affect both the
parabolic Harnack inequality and the heat kernel estimates, see
Theorems \ref{harnackgen2} and \ref{heatgen2}. At the same time it
seems that the Sobolev inequality is independent
 of the $\alpha_k$'s, in the sense that the exponents $\alpha_k$
 do not appear in the ratio (\ref{3.9zero}). We note however that
there are critical cases where relation (\ref{3.9zero}) fails and
different Sobolev inequalities hold true. For instance if
$\alpha_n=-\frac{n-2}{n}$ estimate (\ref{3.9zero}) is no longer
true; instead, the optimal Sobolev inequality involves a logarithmic
correction see inequality (\ref{3.19}) in  Theorem \ref{thm3.2}. For
other examples see Theorem $A'$ in
 \cite{FT}.

We note that instead of the uniform ellipticity condition
(\ref{3.ell}) which we assume throughout this work, our method can
also treat degenerate operators for which the following condition
holds
 \be C_0 w(x)
|\xi|^2 \leq \sum_{i,j=1}^{n} a_{ij}(x) \xi_i \xi_j  \leq C_0^{-1}
w(x)|\xi|^2, \ \ ~~~~ \xi \in \R^n \ ,\ee where $w(x)$ is like a
power of the distance function.

Finally we should mention that this work complements and extends our
previous work \cite{FMT1}. There we studied the cases where the
potential $V(x)$ is either $\frac{(n-2)^2}{4|x|^2}$ for a general
bounded domain $\Omega$ or else $\frac{1}{4 {\rm dist}^2(x,\partial
\Omega)}$ for a convex bounded domain $\Omega$. In the second case
the convexity was used in an essential way. Here, even in these two
cases, we improve our results in the first case by allowing
potentials involving distances to more than one point and in the
second case by removing the convexity assumption (see Section 4).

The article is organized as follows. In Section 2 we establish a
Sobolev type inequality, thus giving the proof of Theorem
\ref{thm3.1zero}, starting from the $L^2$ estimate and using the
 behavior of the generalized eigenfunction $\phi_1$. In
Section 3 we study the associated Cauchy-Dirichlet problem and prove
a parabolic Harnack inequality up to the boundary as well as sharp
two sided estimates on the corresponding heat kernel. In particular
we provide the proof of Theorems \ref{harnackgen2} and
\ref{heatgen2}. Finally, in Section 4 we give some examples of
concrete Schr\"odinger operators with singular potentials for which
an $L^2$ inequality holds true and the  behavior of
$\phi_1$ is known. In all these examples the results of the present
work apply.

\setcounter{equation}{0}
\section{From the $L^2$ estimate to Sobolev type inequalities}

In this section starting from the $L^2$ estimate (\ref{3.5}) we
will prove various Sobolev inequalities involving optimal
exponents. In particular we will prove a weighted logarithmic
Sobolev inequality that will be crucial in establishing the
intrinsic ultracontractivity of the semigroup associated with the
operator ${\mathcal L}$ defined in (\ref{defl}).

For $\xd$ small enough we set
\[
 \xG_k^\xd = \{ x \in \xO,
~~~{\rm s.t.}~~ {\rm dist}(x,\xG_{k} )< \xd \}~~~~{\rm  and}~~~
(\xG_k^\xd)^c= \xO \setminus \xG_k^\xd.
\]

 As a consequence of the assumptions on the domain $\Omega$ we made in the introduction,
 we have that  for  $\xd$ small enough
 $\xG_{k,j}^\xd \cap \xG_{l,i}^\xd = \emptyset$  if $k \neq l$,
or $i \neq j$  and $\xG_{k,j}^\xd \cap \xG_n^\xd=
 \emptyset$ , for $k=1, \ldots, n-1$
and $j=1,\ldots,m_k$.
 We  note that
 $\xO=  \left( \cup_{k=1}^{n} \xG_k^\xd \right) \cup
\left( \cup_{k=1}^{n} \xG_k^\xd \right)^c  =
 \left( \cup_{k=1}^{n} \xG_k^\xd \right) \cup  \left( \cap_{k=1}^{n}
 (\xG_k^\xd)^c \right)$.

\smallskip
We are now ready to give the proof of Theorem \ref{thm3.1zero}.
Indeed we will more generally assume instead of (\ref{3.5}) the
following $L^2$ estimate
\smallskip

\be\la{3.5inf} -\infty < \xl_1 := \inf_{u \in C_0^{\infty}(\xO)}
\frac{ Q[u] } { \int_{\xO} u^{2} dx }, \ee

Then we prove:

\begin{theorem}{\bf (Optimal Sobolev type inequality)} \la{thm3.1}
For  $V \in  L^1_{loc}(\xO)$ we assume that (\ref{3.5inf}) holds
and in addition there exists a  ground state $\phi_1  \in
H^1_{loc}(\xO)  \cap  L^{\infty}_{loc}(\xO)$ satisfying
(\ref{2.3x}) for \be\la{3.8} \xa_k > -\frac{k-2}{2} - (n-k)
\frac{q-2}{2(q+2)} \ ,  \ \ k=1,\ldots,n \ee where $2 <  q \leq
\frac{2n}{n-2}$ if $n\ge 3$ and $q>2$ if $ n=2 .$
 We then  have that for every $\lambda>0$
\be\la{3.9} 0 < C(\xO,\alpha_1,\ldots,\alpha_n, q,\lambda) =
 \inf_{u \in C_0^{\infty}(\xO)} \frac{Q[u]+(
\lambda -\lambda_1) \int_{\Omega} u^2 dx } { \left(
\int_{\xO}d^{\frac{q(n-2)-2n}{2}}
 |u|^{q} dx \right)^{\frac{2}{q}}} \ .
\ee In particular when $n \ge 3$ and $q= \frac{2n}{n-2}$
we have that for every $\lambda>0$, 
$$ 0 < C(\xO,\alpha_1,\ldots,\alpha_n, \lambda) =
 \inf_{u \in C_0^{\infty}(\xO)} \frac{Q[u]+(\lambda-\lambda_1) \int_{\Omega} u^2 dx  }
{ \left( \int_{\xO}
 |u|^{\frac{2n}{n-2}} dx \right)^{\frac{n-2}{n}}}.
$$ if 
$$\xa_k > -\frac{k-2}{2} -
\frac{(n-k)}{2(n-1)},~~~~~~~~k=1,\ldots,n \ . $$
\end{theorem}

In the case $\lambda_1>0$ one can take $\lambda=\lambda_1$, thus
proving Theorem \ref{thm3.1zero}.

\smallskip
\noindent {\em Proof of Theorem \ref{thm3.1}:} It is a consequence
of the following estimate for any $v\in C^\infty_0(\Omega)$
\be\la{3.11} \int_{\xO} \phi_1^2
 \left(\sum_{i,j=1}^{n} a_{ij} v_{x_i} v_{x_j}   +   \lambda   v^2  \right) dx
 \geq C
 \left( \int_{\xO}  \phi_1^q  d^{\frac{q(n-2)-2n}{2}}
 |v|^{q} dx \right)^{\frac{2}{q}},
\ee with $C=C(\alpha_1,\ldots,\alpha_n,\Omega,\lambda)>0$. For
convenience we  write
$v_{x_i}$ instead of $\frac{\partial v}{\partial x_i}$. Let
us accept (\ref{3.11}) and give the proof of the Theorem. Clearly
(\ref{3.11}) is valid not only for smooth functions but also for
functions in the completion of $C_0^{\infty}(\xO)$ under the norm
defined by \be \label{norm} ||v||_{H^1_{\phi_1}}:=\left(\int_{\xO}
\phi_1^2 (v^2+ |\nabla v|^2) dx\right)^{1/2}  \ .\ee
 In particular
we can take $v= \frac{u}{\phi_1}$, with $u \in C_0^{\infty}(\xO)$,
in which case we get \bea\la{3.12}  \sum_{i,j}^{n} \int_{\xO} a_{ij} u_{x_i}
u_{x_j} dx
 -2 \sum_{i,j}^{n} \int_{\xO}  a_{ij} u  u_{x_i}   \frac{ (\phi_1)_{x_j}}{\phi_1}
dx +     \sum_{i,j}^{n}   \int_{\xO} a_{ij}   \frac{ (\phi_{1})_{x_i}
(\phi_{1})_{x_j}}{\phi^2_1} u^2 dx    \nonumber  \\
  + \lambda \int_{\xO}  u^2 dx
\geq  C
 \left( \int_{\xO}   d^{\frac{q(n-2)-2n}{2}}
 |u|^{q} dx \right)^{\frac{2}{q}}.
\eea
On the other hand, by standard approximation arguments, (\ref{1.10b})
is valid also for $\psi = \frac{u^2}{\phi_1}$,  with  $u \in C_0^{\infty}(\xO)$.
   For this choice of the test function,
we get that
\be\la{3.13}
 2 \sum_{i,j}^{n} \int_{\xO}  a_{ij} u  u_{x_i}   \frac{ (\phi_{1})_{x_j}}{\phi_1}
dx -  \sum_{i,j}^{n} \int_{\xO} a_{ij}   \frac{ (\phi_{1})_{x_i}
(\phi_{1})_{x_j}}{\phi^2_1} u^2 dx =   \int_{\xO} (V +\xl_1) u^2 dx.
\ee Combining (\ref{3.12}) and (\ref{3.13})  we conclude that for
any $u \in C_0^{\infty}(\xO)$,
\[
Q[u] +(\lambda-\lambda_1) \int_{\Omega} u^2 dx \geq  C \left(
\int_{\xO}d^{\frac{q(n-2)-2n}{2}}
 |u|^{q} dx \right)^{\frac{2}{q}},
\]
which is   the same as  (\ref{3.9}).

It remains to prove (\ref{3.11}). Because of the ellipticity condition
(\ref{3.ell}), estimate (\ref{3.11})  follows from
\be\la{3.15}
\int_{\xO} \phi_1^2
 \left(|\nabla v|^2    +   v^2  \right) dx
 \geq C
 \left( \int_{\xO}  \phi_1^q  d^{\frac{q(n-2)-2n}{2}}
 |v|^{q} dx \right)^{\frac{2}{q}} \ , \   v \in C_0^{\infty}(\xO) \ . \ee
 In view of (\ref{2.3x}) we may replace $\phi_1$ in (\ref{3.15}) by $d_1^{\xa_1} \ldots
d_n^{\xa_n}$.  Estimate (\ref{3.15}) then is true in $
\cap_{k=1}^{n}
 (\xG_{k}^\xd)^c $. This is
  a   consequence of the standard
Sobolev imbedding of functions
 in $H^1( \cap_{k=1}^{n}
 (\xG_{k}^\xd)^c  )$  in $L^q( \cap_{k=1}^{n}
 (\xG_{k}^\xd)^c  )$
  and the fact that $\xd \leq  d_k(x) \leq D_k :=
 {\rm  sup}_{x \in \xO} d_k(x) < \infty$, $k=1, \ldots, n$.

We therefore need to prove that estimate  (\ref{3.15}) is true when
we replace $\xO$  by   $ \cup_{k=1}^{n} \xG_{k}^\xd $.  As a matter
of fact it is enough  to prove that for any $k =1, \ldots,
n$, 
$$\int_{\xG_{k}^\xd} d_k^{2 \xa_k} ( |\nabla v|^2 +
 v^2   ) dx  \geq C \left( \int_{\xG_{k}^\xd} d_k^{\xb_k q}
 |v|^{q} dx \right)^{\frac{2}{q}}, ~~~~~~  v \in
C_0^{\infty}(\xO), $$ where  $\xb_k = \xa_k-1+ \frac{q-2}{2q}n$.
The validity  of this estimate is given in the next main Lemma.

\begin{lemma}\la{la2.2}
Let $\xO \subset \R^n$, $n \geq 2$, be a bounded domain.
 Suppose  that   $\xG_k  \subseteq \partial \xO$  is
a smooth  boundaryless hypersurface  of codimension $k$,  $k=1,
\ldots, n-1$. When $k=n$ we take  $\xG_n$ to be a point. We also
assume that \be\la{2.11} \xb_k = \xa_k-1+ \frac{q-2}{2q}n \ee
where $2 <  q \leq \frac{2n}{n-2}$ if $n\ge 3$, $q>2$ if $n=2.$
 Then,  for any $k=1, \ldots, n$, there exists a
$C=C(\alpha_k,\Omega,\delta,k)>0$ such that for all $\xd>0$ and all
$ v \in  C_0^{\infty}(\xO)$ there holds \be\la{2.13}
\int_{\xG_{k}^\xd} d_k^{2 \xa_k} ( |\nabla v|^2   + v^2 ) dx \geq C
\|d^{\xb_k}_k v \|^2_{L^{q}(\xG_{k}^\xd)}, \ee provided that
\be\la{2.14} \xa_k \neq -\frac{k-2}{2} - (n-k) \frac{q-2}{2(q+2)}.
\ee
\end{lemma}

\smallskip
\noindent {\em Proof:} Let us fix a $k= 1,\ldots, n$.  We will
initially establish the result for $\xd$ small.  For simplicity we
write $d$ instead of $d_k$.  From Lemma 4.2 \cite{FMaT2} (see also
\cite{FMaT1}), we have that if \be\la{2.17}
   1 < Q \leq \frac{n}{n-1},~~~ b=a-1+ \frac{Q-1}{Q}n, ~~~
{\rm and}~~~ a \neq 1-k, \ee then,  for $\xd$ small there exists
a
   $ C>0$
such that  there holds \be\la{2.19}
 C \| d^b w \|_{L^{Q}(\xG_{k}^\xd)}  \leq
  \int_{\xG_{k}^\xd}
d^a |\nabla w| dx + \int_{\partial \xG_{k}^\xd} d^a  |w| dS_x,
~~~~~~~ w \in C_0^{\infty}(\xO \setminus \xG_{k}). \ee We
apply (\ref{2.19}) to the function $w= |v|^{s}$, $s= \frac{2+q}{2}$.
Also, for $\xa_k$, $\xb_k$, and $q$ as in (\ref{2.11}),
 we set
 $$ Q= qs^{-1},~~~~ b= \xb_k s, ~~~~
 a= b+1 - \frac{Q-1}{Q}n = \frac{\xb_k  q}{2}+ \xa_k.
$$ It is easy to check that $a$, $b$, $Q$ thus defined satisfy
(\ref{2.17}). As far as the condition $a \neq 1-k$ is concerned,
when written in terms of $\xa_k$, $q$, $k$ and $n$, it is
equivalent to
\[
\xa_k \neq -\frac{k-2}{2} - (n-k) \frac{q-2}{2(q+2)},
\]
which is precisely (\ref{2.14}). From (\ref{2.19}) we have
\be\la{2.23} C  \|d^{\xb_k} v \|_{L^q(\xG_{k}^\xd)}^{1+\frac{q}{2}}
 = C \| d^b |v|^s \|_{L^{Q}(\xG_{k}^\xd)}
 \leq     \, s \,
\int_{\xG_{k}^\xd} d^a |v|^{s-1} |\nabla v| dx +
  \int_{\partial \xG_{k}^\xd} d^{a} |v|^s dS_x,
\ee for some positive constant $C$. Using Holder's inequality in the
gradient  term of the right hand side we get \bea\la{2.25}
\int_{\xG_{k}^\xd} d^a |v|^{s-1} |\nabla v| dx &  = &
  \int_{\xG_{k}^\xd} d^{\xa_k} |\nabla v| ~~~
 d^{\frac{\xb_k q}{2}}
 |v|^{\frac{q}{2}} dx   \nonumber \\
& \leq &  \|d^{\xa_k} |\nabla v| \|_{L^2(\xG_{k}^\xd)} ~~~
\|d^{\xb_k} v \|_{L^q(\xG_{k}^\xd)}^{\frac{q}{2}}   \nonumber \\
& \leq & c_{\xe}  \|d^{\xa_k} |\nabla v|
\|_{L^2(\xG_{k}^\xd)}^{1+\frac{q}{2}} + \xe \|d^{\xb_k} v
\|_{L^q(\xG_{k}^\xd)}^{1+\frac{q}{2}}. \eea Hence, from (\ref{2.23})
and (\ref{2.25}) we arrive at \be\la{2.27} (C- \xe s) \|d^{\xb_k} v
\|_{L^q(\xG_{k}^\xd)}^{1+\frac{q}{2}} \leq s c_{\xe} \|d^{\xa_k}
|\nabla v| \|_{L^2(\xG_{k}^\xd)}^{1+\frac{q}{2}}+
  \int_{\partial \xG_{k}^\xd} d^{a} |v|^s dS_x.
\ee To continue we will estimate the trace term in (\ref{2.27}).
Using  Holder's inequality
 we have that \bea\la{2.29}
 \int_{\partial \xG_{k}^\xd} d^{a} |v|^s dS_x  &  =  &  \xd^{a}
\int_{\partial \xG_{k}^\xd} |v|^{\frac{2+q}{2}} dx \leq c  \,
\xd^{a} \left(\int_{\partial \xG_{k}^\xd}
 |v|^{\frac{2(n-1)}{n-2}} dx \right)^{\frac{(n-2)(2+q)}{4(n-1)}}
\nonumber \\
& = & c  \,   \xd^{\frac{(\xb_k-\xa_k) q}{2}}
 \|d^{\xa_k} v\|_{L^{\frac{2(n-1)}{n-2}}
(\partial \xG_{k}^\xd)}^{1+\frac{q}{2}}    \nonumber \\
\eea By the trace imbedding \cite{B}, we have that for $u \in
H^1(\xG_k^{\xd})$ \be\la{1.23} \|
u\|^2_{L^{\frac{2(n-1)}{n-2}}(\partial \Gamma_{k}^\xd)} \leq C(n,k)
\| \nabla u \|^2_{L^2(\Gamma_{k}^\xd)} + M \|u
\|^2_{L^2(\Gamma_{k}^\xd)}, \ee where $M=M(n,\Gamma_{k}^\xd)$.
Applying this to $u =  d^{\alpha_k} v$ we get
 \be\la{1.27} \|
d^{\alpha_k} v\|^2_{L^{\frac{2(n-1)}{n-2}}(\partial \xG_k^{\xd})}
\leq C_2 \int_{\xG_k^{\xd}} d^{2 \alpha_k} v^2 dx  + C_2
 \int_{\xG_k^{\xd}} d^{2\alpha_k}| \nabla  v|^2 dx,
\ee with $C_2=C_2(\alpha_k,\xd, k,n)$.

Indeed after some elementary calculations from (\ref{1.23}) we get
for any $\theta>1$

\bea
\| d^{\alpha_k +\theta} v\|^2_{L^{\frac{2(n-1)}{n-2}}(\partial
\xG_k^{\xd})}
 & \leq &  2C(n,k) (\xa_k+\theta)^2 \int_{\xG_k^{\xd}}
d^{2(\alpha_k + \theta)-2} v^2 dx     \nonumber \\
& &  + 2C(n,k)  \int_{\xG_k^{\xd}}
d^{2(\alpha_k + \theta)}| \nabla v|^2 dx
 +M \int_{\xG_k^{\xd}}
d^{2(\alpha_k + \theta)} v^2 dx  \nonumber \\
 & \leq &  (2C(n,k) (\alpha_k + \theta)^2 \xd^{2 \theta -2}+ M \xd^{2 \theta})
\int_{\xG_k^{\xd}} d^{2 \alpha_k} v^2 dx \nonumber \\
& &  +  2C(n,k) \xd^{2 \theta}
\int_{\xG_k^{\xd}} d^{2\alpha_k}| \nabla  v|^2 dx. \nonumber
\eea
 Whence, \be
\label{traccia}
 \| d^{\alpha_k} v\|^2_{L^{\frac{2(n-1)}{n-2}}(\partial \xG_k^{\xd})}
\leq
 (2C(n,k) (\alpha_k  +  \theta)^2 \xd^{ -2}+ M )
\int_{\xG_k^{\xd}} d^{2 \alpha_k} v^2 dx +   2C(n,k)
\int_{\xG_k^{\xd}} d^{2\alpha_k}| \nabla  v|^2 dx, \ee and therefore
(\ref{1.27}). We note that (\ref{1.27}) is valid even for
nonpositive values of $\alpha_k$.

Combining (\ref{2.27}), (\ref{2.29}), (\ref{1.27}) and then
raising to the power $\frac{4}{2+q}$ we easily conclude
(\ref{2.13})  for $\xd$ small.

The general case follows  by noticing that outside $\xG_{k}^\xd$ for
small $\xd$ the corresponding estimate comes from the standard
Sobolev embeddings and the fact that $\xd \leq  d_k(x) \leq D_k :=
 {\rm  sup}_{x \in \xO} d_k(x) < \infty$, $k=1, \ldots, n$.

This completes the proof of the Lemma as well as of  Theorem  \ref{thm3.1}.

\finedim

\smallskip
We note that when $k=n$  condition  (\ref{2.14}) reads $\xa_n
\neq- \frac{n-2}{2}$.  It turns out that the analogue of the
estimate (\ref{2.13}) in case  $k=n$ and   $\xa_n =  -
\frac{n-2}{2}$, involves logarithmic corrections. More precisely
we have:

\begin{lemma}\la{la2.3}
Let $\xO \subset \R^n$, $n \geq 3$, be a bounded domain
 and   $\xG_n= \{x_0 \}  \subset \partial \xO$  be  a point,  such  that
for some $\xd_0$ small
 $\xG_{n}^{\xd_0}=B_{\xd_0}(x_0) \setminus \{x_0 \} \subset \xO$.
We also assume that 
$$ \xa_n =   - \frac{n-2}{2},~~~~~~~ 2 <  q  \leq \frac{2n}{n-2},
~~~~~~ \xb_n = - \frac{n-2}{2}  -1+ \frac{q-2}{2q}n. $$
 Then,
there exists a $C=C(\Omega,\delta,n)>0$ such that for all $\xd>0$
and all $ v \in C_0^{\infty}(\xO)$ there holds \be\la{2.33}
\int_{\xG_{n}^\xd} d_n^{2 -n} ( |\nabla v|^2   +  v^2   ) dx \geq C
\|d^{\xb_n} X^{\frac12+ \frac{1}{q}} v \|^2_{L^{q}(\xG_{n}^\xd)},
\ee where $X=X(\frac{d_n(x)}{D_n})$,  with
 $X(t) = (1-{\rm  ln}t)^{-1}$,  $0<t \leq 1$, and $ D_n :=
 {\rm  sup}_{x \in \xO} d_n(x) < \infty$.
\end{lemma}

\smallskip
\noindent {\em Proof:} As in the previous Lemma it is enough to
give the proof  for $\xd$ small. We may assume that $x_0=0$,
hence, $d_n(x)=|x|$. Also, for simplicity we suppose that $\xd
=1$. Then we recall the following result for any   $w \in
C_0^{\infty}(B_2)$ the following estimate holds \be\la{2.37}
\int_{B_2} |x|^{2 -n} |\nabla w|^2 dx \geq c \left( \int_{B_2}
|x|^{\beta_n q} X^{1 + \frac{q}{2}} |w|^q  dx
\right)^{\frac{2}{q}} = c
 \|d^{\xb_n} X^{\frac12+ \frac{1}{q}} w \|^2_{L^{q}(B_2)}.
\ee This  is  Lemma 3.2 \cite{FT} in the case $q=\frac{2n}{n-2}$
  and Proposition 6.2 \cite{BFT}
in the case where $2  <   q   \leq \frac{2n}{n-2}$. Given a
function $ v \in C_0^{\infty}(\xO)$, we  extent it  from $B_1$ to
the function $\tilde{v}$ supported in
 $B_2$ so that $$\|  \tilde{v }  \|_{B_2}  \leq C \|  v  \|_{B_1}, $$ holds for a
positive constant depending only on $n$, where we denote by
    $\|  v   \|_{B_1}^2 :=
\int_{B_1} |x|^{2-n} (| \nabla v |^2 +  |v|^2) dx$  (note that
away from the origin $v$ is an $H^1$ function).  We next apply
(\ref{2.37}) to $\tilde{v}$ and the result follows easily. We note
that one cannot take a smaller exponent of the logarithmic term
$X$.

$\hfill \Box$

As a consequence of the above Lemma we have:

\begin{theorem}\la{thm3.2}
For  $V \in  L^1_{loc}(\xO)$ we assume that (\ref{3.5inf}) holds
and in addition there exists a  ground state $\phi_1  \in
H^1_{loc}(\xO)  \cap  L^{\infty}_{loc}(\xO)$ satisfying
(\ref{2.3x}) for $k=1,\ldots,n-1$; $n\ge 3$
$$\xa_k > -\frac{k-2}{2} - (n-k) \frac{q-2}{2(q+2)},~~~~~\xa_n =
-\frac{n-2}{2}, ~~~~2 < q \leq \frac{2n}{n-2}. $$
 We then  have that for every $\lambda>0$
\be\la{3.19} 0 < C(\xO,\alpha_1,\ldots,\alpha_n,q,\lambda) =
 \inf_{u \in C_0^{\infty}(\xO)} \frac{Q[u] +(\lambda-\lambda_1)\int_{\Omega} u^2 dx}
{ \left( \int_{\xO}d^{\frac{q(n-2)-2n}{2}} X^{\frac{q}{2}+1}
 |u|^{q} dx \right)^{\frac{2}{q}}}, \ee
where $X=X(\frac{d_n(x)}{D_n})$,  with
 $X(t) = (1-{\rm  ln}t)^{-1}$,  $0<t \leq 1$ and $D_n:=\sup_{x\in \Omega} d_n(x)<\infty$.
In particular by choosing $q= \frac{2n}{n-2}$ we have that for
every
$\lambda>0$ 
$$ 0 < C(\xO,\alpha_1,\ldots,\alpha_n, \lambda) =
 \inf_{u \in C_0^{\infty}(\xO)} \frac{ Q[u]+(\lambda-\lambda_1)\int_{\Omega} u^2 dx }
{ \left( \int_{\xO}  X^{\frac{2(n-1)}{n-2}}
 |u|^{\frac{2n}{n-2}} dx \right)^{\frac{n-2}{n}}}.
$$
\end{theorem}

\smallskip
\noindent {\em  Proof:} The proof is   quite similar  to  the
proof of Theorem \ref{thm3.1},
  where in the
place of Lemma \ref{la2.2}  one uses Lemma  \ref{la2.3} for $k=n$.
We omit further details.

$\hfill \Box$

Concerning the limit case $q=2$ and  $\xa_k =  -\frac{k-2}{2}$ we
have the following

\begin{theorem}\la{th2.2}
For  $V \in  L^1_{loc}(\xO)$ we assume that (\ref{3.5inf}) holds
and in addition there exists a  ground state $\phi_1  \in
H^1_{loc}(\xO)  \cap  L^{\infty}_{loc}(\xO)$ satisfying
(\ref{2.3x}) for
 $\xa_k \geq - \frac{k-2}{2}$,  $k=1, \ldots, n$; $n\ge 2$.
 We then  have that for every $\lambda>0$
\[
0 < C(\xO,\alpha_1,\ldots,\alpha_n, \lambda) =
 \inf_{u \in C_0^{\infty}(\xO)} \frac{ Q[u]+(\lambda-\lambda_1)\int_{\Omega} u^2 dx } {  \int_{\xO} X^2  \frac{u^{2}}{d^2} dx},
\]
where $X=X(\frac{d(x)}{D})$,  with
 $X(t) = (1-{\rm  ln}t)^{-1}$,  $0<t \leq 1$, and $ D :=
 {\rm  sup}_{x \in \xO} d(x) < \infty$.
\end{theorem}

\smallskip
\noindent {\em  Proof:} The proof is quite similar to the proof of
Theorem \ref{thm3.1}, where in the place of Lemma  \ref{la2.2}
one uses Lemma \ref{la2.7} below. We omit further details.

$\hfill \Box$

\begin{lemma}\la{la2.7}
Let $\xO \subset \R^n$, $n \geq 2$, be a bounded domain.
 Suppose  that   $\xG_k  \subseteq \partial \xO$  is
a smooth  boundaryless hypersurface  of codimension $k$,  $k=1,
\ldots, n-1$. When $k=n$ we take  $\xG_n$ to be a point. Then,  for
any $k=1, \ldots, n$, there exists a   $C>0$ such that for all
$\xd>0$
and all $ v \in  C_0^{\infty}(\xO)$ there holds 
$$\int_{\xG_{k}^\xd}  d_k^{-(k-2)} ( |\nabla v|^2   +   v^2 ) dx
\geq C \int_{\xG_{k}^\xd} X^2  d_k^{ -k } v^2   dx   , $$ where
$X=X(\frac{d_k(x)}{D_k})$,  with
 $X(t) = (1-{\rm  ln}t)^{-1}$,  $0<t \leq 1$, and $ D_k :=
 {\rm  sup}_{x \in \xO} d_k(x) < \infty$.
\end{lemma}

\smallskip
\noindent {\em Proof:} This is proved using similar ideas as in
Lemma \ref{la2.3}. We omit further details.

$\hfill \Box$

As a consequence  of Theorems \ref{thm3.1} and \ref{thm3.2} we
obtain

\begin{theorem}\la{thm3.7} {\bf (Weighted log Sobolev)}
For  $V \in  L^1_{loc}(\xO)$ we assume that (\ref{3.5inf}) holds
and in addition there exists a  ground state $\phi_1  \in
H^1_{loc}(\xO)  \cap  L^{\infty}_{loc}(\xO)$ satisfying
(\ref{2.3x}) for $n\ge 2$  with  \be\la{3.60} \xa_k > -\frac{k-2}{2}
-\frac{n-k}{2(n-1)}, ~~~~~k=1,\ldots,n-1, ~~~~~~
 \xa_n \geq -\frac{n-2}{2}.
\ee
Let
\[
A:= max  \{\xa_1, \xa_2, \ldots, \xa_n, 0 \}.
\]
Then, there exists a positive constant $K$ such that for any $\xe$
positive
 and  for any
 $u \in C_0^{\infty}(\xO)$   there holds
\be\la{3.60a}
\int_{\xO} u^2  {\rm ln}
 \left( \frac{|u|}{\| u \|_2  d_1^{ \xa_1} \ldots d_n^{\xa_n}} \right)
dx \leq \xe Q[u] + \left( K - \frac{n+2A}{4} \,  {\rm ln} \xe
\right) \| u \|_2^2; \ee here $ \| u \|_2^2 =  \int_{\xO} |u|^2
dx$.
\end{theorem}
{\em  Proof:} At first we will show that in each $\xG_{k}^\xd$, the
following estimate holds: \be\la{3.61} \int_{\xG_{k}^\xd} \phi_1^2
v^2  {\rm ln} \left(\frac{|v|}{\|v\|_{\xG_k}}\right) dx \leq \xe
\int_{\xG_{k}^\xd}  \phi_1^2 (|\nabla v|^2 + v^2) dx + \left( K -
\frac{n+2 \xa_k^+}{4}\,  {\rm ln} \xe \right) \|v\|^2_{\xG_k}, \ee
where $\xa_k^+ :=\max \{ \xa_k, 0 \}$  and $\|v\|^2_{\xG_k} :=
\int_{\xG_{k}^\xd} \phi_1^2 v^2 dx$.  To this end, let us assume
first that $w$ is normalized  so that  for $d \mu = \phi_1^2 w^2 dx$
one has $\int_{\xG_{k}^\xd} d \mu = \|w\|^2_{\xG_k} =1$. Then, for
$q>2$,  using Jensen's inequality, we have \be\la{3.63}
\int_{\xG_{k}^\xd} \phi_1^2 w^2 \,  {\rm ln} |w|\,  dx =
\frac{1}{q-2} \int_{\xG_{k}^\xd}  {\rm ln} |w|^{q-2} d \mu \leq
\frac{q}{2(q-2)} {\rm ln} \left( \int_{\xG_{k}^\xd} \phi_1^2 |w|^q
dx \right)^{\frac{2}{q}}. \ee To continue, we will use the estimate
\be\la{3.64} \left( \int_{\xG_{k}^\xd}  \phi_1^2 |w|^q dx
\right)^{\frac{2}{q}} \leq C  \int_{\xG_{k}^\xd} \phi_1^2 (|\nabla
w|^2 +  w^2) dx. \ee In case $k=1, \ldots, n-1$ or  $k=n$ and $\xa_n
> - \frac{n-2}{2}$, (\ref{3.64})
 is a direct consequence of (\ref{2.13}), provided that
$d^{2 \xa_k}_k  \leq c d^{q \xb_k}_k$,  that is, $2 \xa_k  \geq q
\xb_k$. In view of the definition of $\xb_k$ (see (\ref{2.11})),
the requirement $2 \xa_k  \geq q \xb_k$
  is equivalent to
$$q(n-2 + 2 \xa_k) \leq 2 (n +2 \xa_k). $$ We note that when $n\ge
3$ if $\xa_k \leq  0$, then we can choose $q = \frac{2n}{n-2}$,
whereas if  $\xa_k > 0$, then the maximum $q$ one can choose
 is $q = \frac{2 (n +2 \xa_k)}{(n-2 + 2 \xa_k)} < \frac{2n}{n-2}$.
 The same choice of $q$ is feasible when $n=2$ and
  $\alpha_1, \alpha_2> 0$.
 Hence, in any case one takes
\be\la{3.65a} q =  \frac{2 (n +2 \xa_k^+)}{(n-2 + 2 \xa_k^+)}, \ee
On the other hand, in case  $k=n$  and  $\xa_n = - \frac{n-2}{2}$,
estimate  (\ref{3.64})  is a direct consequence of  (\ref{2.33})
with  $q= \frac{2n}{n-2}$ if $n\ge 3$; indeed, in this case one
has $q \xb_n = -n$ and clearly $d^{2-n }_n  \leq c d^{-n}_n
X^{\frac{q}{2}+1}$. In particular, in all cases the choice of $q$
is given by (\ref{3.65a}).

>From (\ref{3.63}) and (\ref{3.64})  we get that
\[
\int_{\xG_{k}^\xd} \phi_1^2 w^2 \,  {\rm ln} |w|\,  dx \leq
\frac{q}{2(q-2)} {\rm ln}   \left( C \int_{\xG_{k}^\xd} \phi_1^2
(|\nabla w|^2 + w^2) dx \right).
\]
Using the fact that $ {\rm ln} \theta \leq \xe \theta - {\rm ln} \xe $ for all
$\theta$, $\xe$ positive  we get that  there exists a $K>0$ such that
  for any $\xe>0$
\be\la{3.66} \int_{\xG_{k}^\xd} \phi_1^2 w^2 \,  {\rm ln} |w|\, dx
\leq \xe \int_{\xG_{k}^\xd}  \phi_1^2 (|\nabla w|^2 + w^2) dx + K -
\frac{q}{2(q-2)}  {\rm ln} \xe. \ee Because of (\ref{3.65a}) we have
that
\[
\frac{q}{2(q-2)} = \frac{n+2 \xa_k^+}{4}.
\]
On the other hand given any $v \in C_0^{\infty}(\xO)  $ we apply
(\ref{3.66}) to  $w=\frac{v}{\|v\|_{\xG_k}}$ to conclude
(\ref{3.61}).

We next consider  a $w \in  C_0^{\infty}(\xO)$ normalized by
 $\|w\|^2_{\xO} :=
\int_{\xO} \phi_1^2 w^2 dx=1$. Applying (\ref{3.61}) to this $w$ we
get \bea \int_{\xG_{k}^\xd} \phi_1^2 w^2  {\rm ln} |w| dx &-&
\int_{\xG_{k}^\xd} \phi_1^2 w^2  {\rm ln}  (\|w\|_{\xG_k}) dx \leq
\xe \int_{\xG_{k}^\xd}  \phi_1^2 (|\nabla w|^2 + w^2) dx
\nonumber \\
&  + & \left( K - \frac{n+2 \xa_k^+}{4}\,  {\rm ln} \xe \right)
\|w\|^2_{\xG_k}.  \nonumber \eea Since $\|w\|_{\xG_k} \leq 1$ and
therefore $ {\rm ln}  (\|w\|_{\xG_k}) \leq 0$, we have in
particular that
\[
\int_{\xG_{k}^\xd} \phi_1^2 w^2  {\rm ln} |w| dx \leq \xe
\int_{\xG_{k}^\xd}  \phi_1^2 (|\nabla w|^2 + w^2) dx + \left( K -
\frac{n+2 \xa_k^+}{4}\,  {\rm ln} \xe \right) \|w\|^2_{\xG_k}.
\]
Summing over all $\xG_{k}^\xd$ we get \be\la{3.67} \int_{\cup
\xG_{k}^\xd} \phi_1^2 w^2  {\rm ln} |w| dx \leq \xe \int_{\cup
\xG_{k}^\xd}  \phi_1^2 (|\nabla w|^2 +  w^2) dx + \left( K -
\frac{n+2 A}{4}\,  {\rm ln} \xe \right) \|w\|^2_{\cup \xG_k}. \ee On
the other hand on $\xO \setminus \cup \xG_{k}^\xd$ we have that
$\phi_1 \sim C$ and using  the standard log Sobolev inequality
 we easily  arrive at
\be\la{3.68} \int_{(\cup \xG_{k}^\xd)^c} \phi_1^2 w^2  {\rm ln} |w|
dx \leq \xe \int_{(\cup \xG_{k}^\xd)^c}  \phi_1^2 (|\nabla w|^2 +
w^2) dx + \left( K - \frac{n}{4}\,  {\rm ln} \xe \right)
\|w\|^2_{(\cup \xG_k)^c}, \ee when $n=2$ (\ref{3.68}) holds true for
any $\nu>2$ in place of $n$. Combining (\ref{3.67}) and (\ref{3.68})
we get that for $\|w\|^2_{\xO} =1$,
 there holds
\be\la{3.69} \int_{\xO} \phi_1^2 w^2  {\rm ln} |w| dx \leq \xe
\int_{\xO}  \phi_1^2 (|\nabla w|^2 + w^2) dx + \left( K -
\frac{n+2 A}{4}\,  {\rm ln} \xe \right). \ee For a general $v \in
C_0^{\infty}(\xO)$ we apply (\ref{3.69}) to $w=
\frac{v}{\|v\|_{\xO}}$ to get \be\la{3.70} \int_{\xO} \phi_1^2 v^2
{\rm ln} \left(\frac{|v|}{\| v \|_{\xO}}\right) dx \leq \xe
\int_{\xO} \phi_1^2 (|\nabla v|^2 + v^2) dx + \left( K - \frac{n+2
A}{4}\,  {\rm ln} \xe \right)\| v \|^2_{\xO} . \ee Taking $u=
\phi_1 v$,   (\ref{3.60a}) follows.

\finedim

The above logarithmic Sobolev inequality is the main ingredient in
establishing the intrinsic ultracontractivity of the semigroup
generated by the operator $\mathcal L$ defined in (\ref{defl}).
More precisely we have

\begin{proposition}\label{intrins}
Let $V\in L^1_{loc}(\Omega)$. We assume that  (\ref{3.5inf}) and
(\ref{2.3x}) hold for some $\xa_k
> -\frac{k-2}{2} -\frac{n-k}{2(n-1)}$, $k=1,\ldots,n-1$, $\xa_n \geq
-\frac{n-2}{2}$. Then the operator ${\mathcal L}$ defined in
(\ref{defl}) gives rise to an intrinsic ultracontractive semigroup
in $L^2(\Omega)$, whose heat kernel $h(t,x,y)$ satisfies \be
\label{up} h(t,x,y)\le C \frac{d_1^{\alpha_1}(x)\cdots
d_n^{\alpha_n}(x) d_1^{\alpha_1}(y) \cdots
d_n^{\alpha_n}(y)}{t^{\frac{n+2A}{2}}} e^{-\lambda_1 t} \ \ \hbox
{for any } t>0, x, y \in \Omega \ ;\ee where $A:= max \{\xa_1,
\xa_2, \ldots, \xa_n, 0 \}.$
\end{proposition}

\noindent {\em Proof:} This is quite similar to Theorem 3.4 in
\cite{FMT1} for this reason we only sketch it. We change variables
by \be \label{change} v(x,t):=u(x,t)/\phi_1(x) , \ee then if $u$
solves problem (\ref{defl}) the function $v$ satisfies \be
\label{deflp} \left\{\begin{array}{lll} \frac{\partial v}{\partial
t}=-{\mathcal L}_{\phi_1} v:=\frac{1}{\phi_1^2}\sum^{n}_{i,j=1}
\frac{\partial}{\partial x_i} \left(\phi_1^2 a_{ij}(x)
\frac{\partial
v}{\partial x_j}\right)- \lambda_1 v   & \hbox { in } \ \ (0,\infty)\times \Omega ~ ,\\
v(x,t)=0 & \ \hbox
{ on } \ \ (0,\infty)\times \partial_1\Omega , \\
v(x,0)=v_0(x)  & \ \hbox { on } \ \ \Omega \ , \end{array} \right.
\ \ \ee with $v_0(x):=u_0(x)\phi_1^{-1}(x)$.

We note that the elliptic operator ${\mathcal L}_{\phi_1}-\lambda_1$
is defined in the domain $D({\mathcal L}_{\phi_1}-\lambda_1):=\{v\in
H^1_0(\Omega;\phi_1^2): {\mathcal L}_{\phi_1}-\lambda_1 \in
L^2(\Omega, \phi_1^2(y) dy)\}$, where $H^1_0(\Omega; \phi_1^2)$
denotes the closure of $C^\infty_0(\Omega)$ functions with respect
to the norm (\ref{norm}). To this elliptic operator it is naturally
associated a bilinear symmetric form which is a Dirichlet form. Then
Lemma 1.3.4 together with Theorems 1.3.2 and 1.3.3 in \cite{D3}
imply that the elliptic operator ${\mathcal L}_{\phi_1}-\lambda_1$
generates
 an analytic semigroup, $e^{-({\mathcal L}_{\phi_1}-\lambda_1) t}$,
 which is
positivity preserving and contractive in $L^p(\Omega,\phi_1^2 dx)$
for any $1\le p\le \infty$.

>From the weighted logarithmic Sobolev inequality (\ref{3.70}), we
deduce the corresponding $L^p$ logarithmic Sobolev inequality for
any $p>2$; to this end it is enough to apply (\ref{3.70})  to the
function $v:=w^\frac{p}{2}$ for any smooth $w$. Using Theorem 2.2.7
in \cite{D3} - as it is used in Corollary 2.2.8  of \cite{D3} - we
obtain that $e^{-({\mathcal L}_{\phi_1}-\lambda_1) t}$ is
 an ultracontractive semigroup. As a consequence the semigroup $e^{-{\mathcal
 L}_{\phi_1}t}$
 has a heat kernel
$h_{\phi_1}$, which satisfies the following uniform upper bound
\be \label{uppp} h_{\phi_1}(t,x,y)\le
\frac{C}{t^{\frac{n+2A}{2}}}e^{-\lambda_1 t} \ \ \hbox {for any }
t>0, x, y \in \Omega \ .\ee Clearly the heat kernel upper bound
(\ref{up}) associated to the operator $\mathcal L$ follows from
(\ref{uppp}), (\ref{2.3x}) and the fact that \be\label{changebis}
h(t,x,y)=\phi_1(x)\phi_1(y) h_{\phi_1}(t,x,y) \ ,\ee which is an
immediate consequence of the change of variables (\ref{change}).
We omit further details.

\finedim

\section{Harnack inequalities and sharp heat kernel estimates}
\setcounter{equation}{0}

In this section we prove a parabolic Harnack inequality up to the
boundary for the operator ${\mathcal L}_{\phi_1}$ defined in
(\ref{deflp}), and we deduce from it the corresponding heat kernel
estimates as well as the proofs of Theorems \ref{harnackgen2} and
\ref{heatgen2} in the Introduction. We use Moser iteration
technique, as adapted in \cite{FMT1} for bounded domains $\Omega$.
To this end we will prove four basic estimates. Namely, a sharp
volume estimate, a local weighted Poincar\'e inequality, a local
weighted Moser inequality and a density theorem.

We will use the following local representation of any smooth
boundaryless hypersurface $\xG_{k,j}$ of codimension
$k=1,\ldots,n-1$, for any fixed $j=1, \cdots, m_k$, which is of
course Lipschitz. That is, we suppose there exists a finite number
$N$ (depending on both $k$ and $j$) of coordinate systems
$(y_i,z_i)$, $y_i=(y_{i,1}, \cdots, y_{i,(n-k)})$ and
$z_i=(z_{i,1}, \cdots, z_{i,k})$, for $i=1,\cdots , N$, and the
same number of functions $a_i=a_i(y_i):\R^{n-k}\to \R^k$,
($a_i=(a_i^1,\ldots, a_i^k)$) defined on the closures of the
$(n-k)$ dimensional cubes $\Delta_i:=\{y_i:|y_{i,l}|\le \beta
\hbox { for } l=1,\cdots, n-k\}$, $i\in \{1,\cdots , N\}$ so that
for each point $x\in \xG_{k,j}$ there is at least one $i$ such
that $x=(y_i, a_i(y_i))$. The functions $a_i$ satisfy the
Lipschitz condition on $\overline\Delta_i$ with a constant $L>0$
that is
$$|a_i(y_i)-a_i(\bar y_i)|_{\R^k}\le L|y_i-\bar y_i|_{\R^{n-k}}$$ for any
$y_i, \bar y_i \in \overline \Delta_i$. We note $|y|_{\R^{k}}$ is
the Euclidean norm in $\R^k$. Moreover, there exists a positive
number $\beta<1$, called the localization constant of $\Gamma_{k,j}$
and $\Omega$, such that the set $B_i$ defined for any $i\in\{1,
\cdots, N\}$ by the relation
$$B_i=\{(y_i,z_i): y_i\in \Delta_i, \ \ a_i^l(y_i)-\beta
<z_{i,l}<a_i^l(y_i)+\beta\} \ , $$ satisfies
$$U_i=B_i \cap
\Omega = \{(y_i,z_i): y_i\in \Delta_i, \ \ a_i^l(y_i)-\beta
<z_{i,l}<a_i^l(y_i)\} \hbox { if } k=1 \ ,$$ or
$$U_i=B_i \cap \Omega=B_i =
\{(y_i,z_i): y_i\in \Delta_i, \ \ a_i^l(y_i)-\beta
<z_{i,l}<a_i^l(y_i)+\beta\} \ \hbox{ if } k= 2,\ldots,n-1 \ ,$$
and $\Gamma_i=B_i \cap
\partial \Omega=\{(y_i,z_i): y_i\in \Delta_i,
z_i=a_i(y_i)\}$. Finally, let us observe that for any $y\in U_i$
one has $$(1+L)^{-1} |a_i(y_i)-z_i|_{\R^k} \le d_k(y_i,z_i)\le
|a_i(y_i)-z_i|_{\R^k} \ ;$$ see Corollary 4.8 in \cite{K}.

\smallskip
\smallskip

We fix a constant $\gamma \in (1,2)$ and we define the ``balls" we
will use in Moser iteration technique. Roughly speaking they will
be Euclidean balls if they stay away from the boundary and they
will be $n$ dimensional ``deformed cubes", following the geometry
of the boundary, if they are close enough to the boundary or if
they intersect it. More precisely we have

\begin{definition}\label{ball}
(i) For any $x\in \Omega$ and for any $0<r<\beta$ we define the
``ball" centered at $x$ and having radius $r$ as follows.
${\mathcal B}(x,r)=B(x,r)$ the Euclidean ball centered at $x$ and
having radius $r$ if $d(x)\ge \gamma r$ (thus $d_k(x)\ge \gamma r$
for any $k=1, \cdots, n$) or if $k=n$, while
$${\mathcal B}(x,r)=\{(y_i,z_i): |y_{i}-x'|_{\R^{n-k}} < r ,
\ $$ $$ a_i^l(y_i)-r-d_k(x) <z_{i,l}<a_i^l(y_i)+r-d_k(x) \hbox {
for any } l=1,\cdots, k \}$$ if $k=1,\ldots, n-1$ and
$d_k(x)<\gamma r$ where $i\in \{1, \cdots , N\}$ is uniquely
defined by the point $\bar x \in \xG_k$ such that $|\bar
x-x|_{\R^n}=d_k(x)$, that is by the projection of the center $x$
onto $\xG_k\subset
\partial \Omega$, and  $x'$ denotes the first $n-k$ coordinates of the point $x$ in the
$i$-orthonormal coordinate system. (ii) We also define the volume
of the ``ball" centered at $x$ and having radius $r$ by
$$V(x,r):=\int_{{\mathcal B}(x,r)\cap \Omega} \prod^n_{k=1}
d_k^{2\alpha_k}(y) dy \ .$$
\end{definition}

\bigskip

We first derive a sharp volume estimate, which plays a fundamental
role in getting the sharp dependence of the heat kernel on $x,y$
and $t$.

\begin{lemma}
\label{volume}{\bf (Sharp volume estimate)} Let $n\ge 2$
 and $\alpha_k>-\frac{k}{2}$ for $k=1,\ldots,n$. Then, there exist positive constants
$c_1, c_2$ and $r_0$ such that for any $x\in \Omega$ and
$0<r<r_0$, we have \be\label{volums} c_1 \prod^n_{k=1} (d_k (x)+
r)^{2\alpha_k} r^n \le V(x,r)\le c_2 \prod^n_{k=1} (d_k (x)+
r)^{2\alpha_k} r^n \ .\ee
\end{lemma}

\smallskip
\noindent {\em Proof of Lemma \ref{volume}:} Let us first consider
the case where $d(x)\ge \gamma r$, whence $d_k(x)\ge \gamma r$ for
any $k=1,\ldots,n$. Then ${\mathcal B}(x,r)=B(x,r)\subset \Omega$.
Due to the fact that for any $y\in B(x,r)$ and any $k=1,\cdots, n$
we have \be\label{bee} \left(\frac{\gamma-1}{\gamma}\right)
d_k(x)\le d_k(x)-r\le d_k(y)\le d_k(x)+r \le
\left(\frac{\gamma+1}{\gamma}\right) d_k(x) \ee we easily get
$$V(x,r)\sim r^n \prod^n_{k=1} d_k^{2\alpha_k}(x) \sim r^n \prod^n_{k=1} (d_k (x)+ r)^{2\alpha_k}\ ,$$
this proves the claim.

Let us now consider the case where $d(x)<\gamma r$. We claim that in
this case there exists exactly one $k=1,\cdots, n$ such that
$d_k(x)<\gamma r$. This is due to the assumption that for some
$\delta$ small enough $\xG_{k,j}^\delta\cap
\xG_{l,i}^\delta=\emptyset$ for any $k\neq l$ and $i\neq j$ and
$\xG_{k,j}^\delta\cap \xG_{n}^\delta=\emptyset$ for any $k=1,\cdots,
n-1$ and $j=1,\cdots, m_k$, since we may suppose that
$r<\frac{\delta}{2}$ (take $r_0:=\min\{\beta,\frac{\delta}{2}\}$,
$\beta$ being the localization constant of $\Gamma_{k,j}$ and
$\Omega$). Whence if $d_k(x)<\gamma r$ then $x\in \xG_{k}^\delta$
and $d_j(x) \ge \delta
> \gamma r$ for any $j\neq k$, thus from (\ref{bee}) for any $y\in
B(x,r)$ we have $d_j(y)\sim d_j(x)$; as a consequence
$$V(x,r)\sim \prod_{{j=1, j\neq k}}^n (d_j(x)+r)^{2\alpha_j}
\int_{{\mathcal B}(x,r) \cap \Omega} d_k^{2\alpha_k}(y) dy \ ,$$
Hence the claim will follow as soon as we prove that \be
\label{unterm}\int_{{\mathcal B}(x,r) \cap \Omega}
d_k^{2\alpha_k}(y) dy \sim r^{n+2\alpha_k} \ .\ee Arguing as in
(\ref{bee}) we have that $d_k(y) \le (1+\gamma) r$ for any $y\in
{\mathcal B}(x,r)$. Moreover one has that if $k\neq n$, $d_k(y)\ge
r(\gamma-1)$ on a set of measure $r^n$. Indeed
$$\int_{{\mathcal B}(x,r)\cap \Omega \cap \{d_k(y)\ge
r(\gamma-1)\}} dy = \int_{|y_{i}-x'|_{\R^k}\le r}
\int^{a_i^l(y_i)+r-r\gamma}_{a_i^l(y_i)-r-d_k(x)} dz_{i,l} dy_i
=$$
$$=(2r-\gamma r+d_k(x))^k r^{n-k}\ge (2-\gamma)^k r^{n} ,$$ and (\ref{unterm}) follows.
In the limit case $k=n$ see \cite{MT2} for any $\alpha_n \in
(-\frac{n}{2},0]$ (see also \cite{MT1} and Lemma 2.3 in
\cite{FMT1}). We note in fact that the same proof works for any
$\alpha_n>-\frac{n}{2}$.

\finedim

>From this one can easily deduce the doubling property:

\begin{corollary}{\bf (Doubling property)}
Let $n\ge 2$ and $\alpha_k>-\frac{k}{2}$ for $k=1,\ldots,n$. Then,
there exist positive constants $C_D$ and $r_0$ such that for any
$x\in \Omega$ and $0<r<r_0$, we have
$$V(x,2r)\le C_D V(x,r) \ .$$
\end{corollary}

\bigskip

Our next result reads:

\begin{theorem}{\bf (Local weighted Poincar\'e inequality)}
\label{poincare} Let  $n\ge 2$ and $\alpha_1> 0$, $\alpha_k >
-\frac{k}{2}$ for $k=2,\cdots,n$. Then, there exist positive
constants $C_P$ and $r_0$ such that for any $x\in \Omega$ and
$0<r<r_0$, we have for all $f \in C^1(\overline {{\mathcal
B}(x,r)\cap \Omega})$
\begin{equation} \label{poipoi}\inf_{\xi \in \R}\int_{{\mathcal
B}(x,r)\cap \Omega} \prod_{k=1}^n d_k^{2\alpha_k}(y)|f(y)-\xi|^2
dy\le C_P \ r^2 \int_{{\mathcal B}(x,r)\cap \Omega} \prod_{k=1}^n
d_k^{2\alpha_k}(y) |\nabla f|^2 dy \ .
\end{equation}
\end{theorem}

We note that our weight is not necessarily in the Muckenhoupt
class $A^2$.

\bigskip
\noindent {\em Proof: } Let us first consider the case where
$d(x)\ge \gamma r$. Then ${\mathcal B}(x,r)=B(x,r)\subset\Omega$
and for any $y\in B(x,r)$ and any $k=1,\cdots, n$ we have
$d_k(y)\sim d_k(x)$, as in estimate (\ref{bee}). Thus in this case
(\ref{poipoi}) follows from the standard Poincar\'e inequality:
$$\inf_{\xi \in \R}\int_{B(x,r)}
|f(y)-\xi|^2  dy\le C_P \ r^2 \int_{B(x,r)} |\nabla f|^2 dy \ , ~
~ ~~~~ \ f \in C^1(\overline {B(x,r)}) \ .
$$

Let us now consider the case where $d_k(x)<\gamma r$, for some
$k=1,\cdots, n$. Then arguing as in Lemma \ref{volume} it is enough
to prove the following for any $f \in C^1(\overline {{\mathcal
B}(x,r)\cap \Omega})$ and any $k=1,\cdots, n$ \be
\label{epoi}\inf_{\xi \in \R}\int_{{\mathcal B}(x,r)\cap \Omega}
|f(y)-\xi|^2 d_k^{2\alpha_k}(y)  dy\le C_P \ r^2 \int_{{\mathcal
B}(x,r)\cap \Omega} |\nabla f|^2 d_k^{2\alpha_k}(y) dy \ . \ee The
case $k=1$ corresponds to Theorem 2.5 in \cite{FMT1} (with
$\lambda=0$ there). The case $k=n$ has been treated in Theorem 3.1
in \cite{MT2} (see also \cite{MT1} and Theorem 2.5 in \cite{FMT1})
  for any $\alpha_n\in (-\frac{n}{2},0]$. We note however that the same proof works
  for any $\alpha_n>-\frac{n}{2}$. So we need to consider the intermediate cases $k=2, \ldots,
  n-1$.

\smallskip
\noindent  We deduce (\ref{epoi}) from the analogous statement for
$k=n$, that is from the following inequality
 \be
\label{epoibo}\inf_{\xi \in \R}\int_{B(x,r)} |f(y)-\xi|^2
|y|_{\R^n}^{2\alpha_n} dy\le C_P \ r^2 \int_{B(x,r)} |\nabla f|^2
|y|_{\R^n}^{2\alpha_n} dy, \ \ ~~~ \ f \in C^1(\overline
{B(x,r)}) . \ee As a consequence of the local representation we have
for some $a$ and $s=(y,z)$
$$\int_{{\mathcal B}(x,r)\cap \Omega}d_k^{2\alpha_k}
 (s)|f(s)-\tilde \xi|^2  ds \le $$ $$\le C(L)\int_{|y-x'|_{\R^{n-k}}\le r} \int^{
a^l(y)+r-d_k(x)}_{a^l(y)-r-d_k(x)} |f(y,z)-\tilde \xi|^2
|a(y)-z|_{\R^k}^{2\alpha_k} dz_l dy \le$$
$$\le C \int_{|y-x'|_{\R^{n-k}}\le r}\int_{d_k(x)-r}^{r+d_k(x)} |f(y,a(y)-w)-\tilde \xi|^2
|w|_{\R^k}^{2\alpha_k} \ dw_l dy \ ,$$ here we used the following
change of variables $(y,z)\to (y,w:=a(y)-z)$. Then, since
$$|f-\tilde \xi|^2\le 2\left( |f-\xi(y)|^2+|\xi(y)-\tilde \xi|^2\right) \ ,$$ where we use the following
notation
$$\xi(y):=\frac{\int_{d_k(x)-r}^{r+d_k(x)} f(y,a(y)-w)
|w|_{\R^k}^{2\alpha_k} dw }{\int_{d_k(x)-r}^{r+d_k(x)}
|w|_{\R^k}^{2\alpha_k} dw } \ ,$$
$$\tilde \xi:=\omega_{n-k-1}^{-1}r^{-n+k}\int_{|y-x'|_{\R^{n-k}}\le
r} \xi(y) dy \ ,$$ inequality (\ref{epoi}) follows from estimates
(i) and (ii) below.

(i) We have
$$\int_{|y-x'|_{\R^{n-k}}\le r} \left(\int_{d_k(x)-r}^{r+d_k(x)} |f(y,a(y)-w)-\xi(y)|^2
|w|_{\R^k}^{2\alpha_k} \ dw_l \right) dy \le $$ $$\le C
(r+d_k(x))^2 \int_{|y-x'|_{\R^{n-k}}\le r}
\left(\int_{d_k(x)-r}^{r+d_k(x)} |\nabla_z f|^2
|w|_{\R^k}^{2\alpha_k} \ dw_l \right) dy \le $$ $$\le C r^2
\int_{{\mathcal B}(x,r)\cap \Omega} d_k^{2\alpha_k}(s)|\nabla f|^2
ds \ ,$$ here we used the assumption $d_k(x)<\gamma r$ as well as
inequality (\ref{epoibo}) applied in $\R^k$, instead of $\R^n$,
this explains the restriction $2\alpha_k > -k$.

(ii) Finally
$$\int_{d_k(x)-r}^{r+d_k(x)} \left(\int_{|y-x'|_{\R^{n-k}}\le r} |\xi(y)-\tilde \xi|^2
dy \right)|w|_{\R^k}^{2\alpha_k} \ dw_l  \le $$ $$\le C r^2
\int_{d_k(x)-r}^{r+d_k(x)} \int_{|y-x'|_{\R^{n-k}}\le r}
\Big|\nabla_y f+\sum^{k}_{l=1} \frac{\partial f}{\partial z_l}
\nabla_y a^l(y)\Big|^2 dy |w|_{\R^k}^{2\alpha_k} \ dw_l \le$$
$$\le
C r^2 \int_{{\mathcal B}(x,r)\cap \Omega} d_k^{2\alpha_k}|\nabla
f|^2 dydz \ ,$$ here we used the standard Poincar\'e inequality on
the Euclidean $n-k$ dimensional ball of radius $r$ centered at
$x'$.

\smallskip
\noindent The proof of inequality (\ref{epoi}) is now complete.

\finedim

All  the ingredients of the abstract machinery of  \cite{GSC} are
now in place. However, since  bounded domains endowed with the
Euclidean metric are not complete manifolds, the standard method
should be modified as in \cite{FMT1}. In particular we will next
prove a local weighted Moser inequality as well as  a density
Theorem which are crucial in making the Moser iteration to work in
our setting.

\bigskip

We next prove the following local weighted Moser inequality:

\smallskip
\begin{theorem}
{\bf (Local weighted Moser inequality)} \label{moser} Let $n\ge 2$
 and $\alpha_1>0 \ , \alpha_k \ge -\frac{k-2}{2}$  for $k=2,\cdots, n$.
 Then, there exist positive constants $C_M$ and
$r_0$ such that for any $\nu\ge n+2A$, $A:= max \{\xa_1, \xa_2,
\ldots, \xa_n, 0 \}$, $x\in \Omega$, $0<r<r_0$ and $f \in
C^\infty_0({\mathcal B}(x,r)\cap \Omega)$ we have \be
\label{momo}\int_{{\mathcal B}(x,r)\cap
\Omega}|f(y)|^{2\left(1+\frac{2}{\nu}\right)} \prod^n_{i=1}
d_i^{2\alpha_i}(y)  dy \le $$ $$ \le C_M r^2
V(x,r)^{-\frac{2}{\nu}} \left(\int_{{\mathcal B}(x,r)\cap
\Omega}|\nabla f|^2 \prod^n_{i=1} d_i^{2\alpha_i}(y)dy \right)
\left(\int_{{\mathcal B}(x,r)\cap \Omega} f^2 \prod^n_{i=1}
d_i^{2\alpha_i}(y) dy \right)^\frac{2}{\nu} . \ee
\end{theorem}

\smallskip
\noindent {\em Proof: } Let us first consider the case where
$d(x)\ge \gamma r$. Then ${\mathcal B}(x,r)=B(x,r)\subset\Omega$ and
for any $y\in B(x,r)$ and any $k=1,\cdots, n$ we have $d_k(y)\sim
d_k(x)$, as in estimate (\ref{bee}), and the claim follows from the
standard Moser inequality, which we recall here: There exists a
positive constant $C$ such that for any $x\in \Omega$, $r>0$, and
any $\nu\ge n$ if $n\ge 3$ or any $\nu>2$ if $n=2$, the following
holds true
$$\int_{B(x,r)}|f(y)|^{2\left(1+\frac{2}{\nu}\right)} dy\le C r^2
r^{-\frac{2n}{\nu}} \left(\int_{B(x,r)} |\nabla f|^2 dy \right)
\left(\int_{B(x,r)} f^2 dy \right)^\frac{2}{\nu} \ \ , $$ for all
$f \in C^\infty_0(B(x,r))$ (see for example Section 2.1.3 in
\cite{SC}). Making use of the sharp volume estimate in Lemma
\ref{volume}, we have
$$\int_{B(x,r)}\prod^n_{i=1} d_i^{2\alpha_i}(y) |f(y)|^{2\left(1+\frac{2}{\nu}\right)}  dy
\le C \prod^n_{i=1} d_i^{2\alpha_i}(x) r^2 r^{-\frac{2n}{\nu}}
\left(\int_{B(x,r)} |\nabla f|^2 dy \right) \left(\int_{B(x,r)}
f^2 dy \right)^\frac{2}{\nu}\le
$$
$$\le C \left(\prod^n_{i=1} d_i^{2\alpha_i}(x)\right)^{1-1-\frac{2}{\nu}} r^2 r^{-\frac{2n}{\nu}}
\left(\int_{B(x,r)} \prod^n_{i=1} d_i^{2\alpha_i}(y) |\nabla f|^2
dy \right) \left(\int_{B(x,r)} \prod^n_{i=1} d_i^{2\alpha_i}(y)
f^2 dy \right)^\frac{2}{\nu} =
$$
$$= C_M r^2 V(x,r)^{-\frac{2}{\nu}} \left(\int_{B(x,r)} \prod^n_{i=1} d_i^{2\alpha_i}(y)
 |\nabla f|^2 dy \right)
\left(\int_{B(x,r)} \prod^n_{i=1} d_i^ {2\alpha_i}(y) f^2 dy
\right)^\frac{2}{\nu}\ ,$$ and (\ref{momo}) has been proved in
case $d(x)\ge \gamma r$.

Let us now consider the case where $d(x)<\gamma r$. Arguing as in
Lemma \ref{volume} this corresponds to consider the case where
$d_k(x)<\gamma r$ for some $k=1,\ldots,n$. In view of
(\ref{volums}) it is enough to prove  \be \label{interest}
\int_{{\mathcal B}(x,r)\cap
\Omega}|f(y)|^{2\left(1+\frac{2}{\nu}\right)} d_k^{2\alpha_k}(y)
dy \le $$ $$ \le C_M r^2 r^{-\frac{2(n+2\alpha_k)}{\nu}}
\left(\int_{{\mathcal B}(x,r)\cap \Omega}|\nabla f|^2
d_k^{2\alpha_k}(y)dy \right) \left(\int_{{\mathcal B}(x,r)\cap
\Omega} f^2 d_k^{2\alpha_k}(y) dy \right)^\frac{2}{\nu} . \ee In
the argument that follows we omit the integral set which is always
taken as ${\mathcal B}(x,r)\cap \Omega$ and we define
$d\mu:=d_k^{2\alpha_k}(y) dy$. First of all making use twice of
H\"older inequality for any $\nu> n+2\alpha_k^+$, we have
\be\label{iniz}\int f^{2\left(1+\frac{2}{\nu}\right)} d\mu\le
\left(\int f^{2\left(1+\frac{2}{n+2\alpha_k^+}\right)}
d\mu\right)^{\frac{\nu+2}{\nu}\frac{n+2\alpha_k^{+}}{n+2\alpha_k^{+}+2}}
\left(\int d\mu
\right)^{1-\frac{\nu+2}{\nu}\frac{n+2\alpha_k^{+}}{n+2\alpha_k^{+}+2}}
\ , \ee as well as \be \label{iniz2} \left(\int f^{2} d\mu
\right)^{\frac{2(\nu-n-2\alpha_k^{+})}{\nu(n+2\alpha_k^+)}}\le
\left(\int f^{2\left(1+\frac{2}{n+2\alpha_k^+}\right)} d\mu
\right)^{\frac{2(\nu-n-2\alpha_{k}^{+})}{\nu(n+2\alpha_k^{+}+2)}}
\left(\int d\mu
\right)^{\frac{4(\nu-n-2\alpha_k^{+})}{\nu(n+2\alpha_k^{+})(n+2\alpha_k^{+}+2)}}\
,\ee multiplying both sides of inequalities (\ref{iniz}) and
(\ref{iniz2}) we deduce that \be \label{tom} \int
f^{2\left(1+\frac{2}{\nu}\right)} d\mu \le \left(\int
 f^{2\left(1+\frac{2}{n+2\alpha_k^+}\right)}
d\mu\right)
 \left(\int f^{2}
d\mu \right)^{\frac{2(n+2\alpha_k^{+}-\nu)}{\nu(n+2\alpha_k^+)}}
\left(\int d\mu
\right)^{\frac{2(\nu-n-2\alpha_k^{+})}{\nu(n+2\alpha_k^{+})}}
 \ .\ee
H\"older inequality also implies that \be \label{mer} \int
f^{2\left(1+\frac{2}{n+2\alpha_k^+}\right)} d\mu \le \left(\int
f^{\frac{2(n+2\alpha_k^+)}{n+2\alpha_k^+-2}}
d\mu\right)^{\frac{n+2\alpha_k^+-2}{n+2\alpha_k^+}}\left(\int f^2
d\mu \right)^{\frac{2}{n+2\alpha_k^+}} \ . \ee To continue we will
use the following local weighted Sobolev inequality
\begin{equation} \label{ls} \left(\int_{{\mathcal B}(x,r)\cap
\Omega} d_k^{2\alpha_k}(y)
|f(y)|^{\frac{2(n+2\alpha_k^+)}{n+2\alpha_k^+-2}}
dy\right)^{\frac{n+2\alpha_k^+-2}{n+2\alpha_k^+}}\le C_S
\int_{{\mathcal B}(x,r)\cap \Omega}  d_k^{2\alpha_k}(y)|\nabla
f|^2 dy \ .\end{equation} Then from  (\ref{tom}), (\ref{mer}),
(\ref{ls}) and (\ref{unterm}) we get the desired Moser inequality
(\ref{interest}) with
$$C_M=C_S \ r^{-\frac{4\alpha^-_{k}}{n+2\alpha_k^+}} \ ,$$ where
$\alpha_k^{-}:=\max\{0, -\alpha_k\}$. It remains to prove
(\ref{ls}) with a positive constant $C_S$ independent of $x$ and
$r$ for $0<r<r_0$. This is a consequence of (\ref{3.64}) and the
local weighted Poincar\'e inequality (\ref{poipoi}) (since for
functions $f\in C^\infty_0({\mathcal B}(x,r)\cap \Omega)$ one can
take $\xi=0$ in it). It follows that $C_M$ is automatically
independent of $r$ if $\alpha_k\ge 0$. Hence Theorem \ref{moser}
has been proved in this case.

It remains to show that $C_M$ is independent of $r$ also in the
case $\alpha_k<0$. In fact in such a case instead of (\ref{ls}) we
have an even better estimate (see the definition of $\beta_k$
given in (\ref{2.11}) for $q=\frac{2n}{n-2}$ and $n\ge 3$)
$$\left(\int_{{\mathcal B}(x,r)\cap \Omega}
d_k^{\frac{2n\alpha_k}{n-2}}(y) |f(y)|^{\frac{2n}{n-2}}
dy\right)^{\frac{n-2}{n}}\le \tilde{C}_S \int_{{\mathcal
B}(x,r)\cap \Omega}  d_k^{2\alpha_k}(y)|\nabla f|^2 dy \ ,$$ for
some positive constant $r_0$ and for any $x\in \Omega$, $0<r<r_0$,
$f\in C^\infty_0({\mathcal B}(x,r)\cap \Omega)$, with a positive
constant $\tilde C_S$ independent from $x$ and $r$. Whence
$$\left(\int_{{\mathcal B}(x,r)\cap \Omega} d_k^{2\alpha_k}(y)
|f(y)|^{\frac{2n}{n-2}} dy\right)^{\frac{n-2}{n}}\le \tilde{C}_S
\left(\max_{y\in {{\mathcal B}(x,r)\cap \Omega}}
d_k(y)\right)^{\frac{4\alpha_k^{-}}{n}} \int_{{\mathcal
B}(x,r)\cap \Omega} d_k^{2\alpha_k}(y)|\nabla f|^2 dy \le$$ $$\le
\tilde{C}_S \ r^{\frac{4\alpha_k^{-}}{n}} \int_{{\mathcal
B}(x,r)\cap \Omega} d_k^{2\alpha_k}(y)|\nabla f|^2 dy \ ,$$ since
$d_k(y)\le d_k(x)+r\le (1+\gamma)r$. If we use the above
inequality in the place of (\ref{ls})  we get $C_M=\tilde C_S$,
that is $C_M$ is independent of $r$.

\smallskip
\noindent The proof of Theorem \ref{moser} is now complete.

\finedim

Finally let us prove the following density Theorem, which as
explained in \cite{FMT1} is crucial for the Moser iteration to work
on bounded domains. Let $$H^1(\xO;\, \phi^2_1)=\{v \in L^2(\Omega,
\phi_1^2 dx): \int_{\xO} \phi_1^2
 |\nabla v|^2 dx < +\infty \}, $$ with norm defined by $
||v||^2_{H^1_{\phi_1}}:=\int_{\xO} \phi_1^2 (v^2+ |\nabla v|^2)
dx$.

\begin{theorem}{\bf (Density theorem)} \la{thm2.6}
Let $n\ge 2$. Suppose that   $\phi_1$  satisfies
\[
c_1 d_1^{\xa_1}(x) \ldots  d_n^{\xa_n}(x)   \leq   \phi_1(x)  \leq
c_2 d_1^{\xa_1}(x) \ldots  d_n^{\xa_n}(x),
\]
for $x\in \Omega$ with   $c_1$, $c_2$    positive constants  and
 $\xa_k \geq - \frac{k-2}{2}$,  $k=1, \ldots, n$.
Then, the
\[
 C^{\infty}_0(\xO)  ~functions ~are ~dense ~in~
  H^1(\xO;\, \phi^2_1).
\]
\end{theorem}

\smallskip
\noindent {\em Proof:} The special case $\alpha_k=0$ for
$k=2,\ldots,n$ was treated in Theorem 2.11 of \cite{FMT1}. First of
all from Theorem 7.2 in \cite{K} it is known that the set
$C^\infty(\overline \Omega)$ is dense in $H^1(\Omega; \phi_1^2)$.
Thus for any $v\in H^1(\Omega; \phi_1^2)$ there exists $v_m\in
C^\infty(\overline \Omega)$ such that for any $\epsilon>0$ we have
$||v-v_m||_{H^1_{\phi_1}}\le \epsilon$ if $m\ge m(\epsilon)$. Let us
choose $w:=v_{m(\epsilon)}$ and let us define, for any $k=1, \cdots,
n$ and any $j\ge 1$, the following function
$$\psi_k^j(x)=\left\{\begin{array}{lll} 0 & \hbox { if } \  d_k(x)\le \frac{1}{j^2}  ~ ,\\
1 +\frac{\ln (jd_k(x))}{\ln (j)} & \hbox { if } \ \frac{1}{j^2}<
d_k(x)<
\frac{1}{j}~ , \\
 1 & \hbox { if }\ d_k(x)\ge \frac{1}{j}  \ .\end{array} \right. $$
Then $w^j:=w \prod^n_{k=1} \psi_k^j\in C^{0,1}_0(\Omega)$, and
$$||w-w^j||_{H^1_{\phi_1}}= ||w
(1-\prod^n_{k=1}\psi_k^j)||_{H^1_{\phi_1}}\le 2\int_{\Omega}
(w^2+|\nabla w|^2) (1-\prod^n_{k=1}\psi_k^j)^2 \phi_1^2(y) \ dy +
$$ $$+2 \int_{\Omega} w^2 \sum^n_{k=1}|\nabla \psi_k^j|^2
\phi^2_1(y) \ dy \le$$
$$\le 2\int_{\cup_{k=1}^n \{d_k(y)< \frac{1}{j}\}} (w^2+|\nabla w|^2) \phi_1^2 dy + 2
\sum^n_{k=1}\int_{\frac{1}{j^2}<d_k(y)< \frac{1}{j}}
\frac{w^2|\nabla d_k|^2}{d_k^2(y) (\ln (j))^2} d_k^{2\alpha_k}(y) dy
\ .$$ Now as $j\to \infty$ it is clear that the first term in the
right hand side goes to zero since $w\in H^1_{\phi_1}$. We next show
that also the second term goes to zero. Recalling that $|\nabla
d_k|\le 1$
$$\int_{\frac{1}{j^2}<d_k(y)< \frac{1}{j}}
\frac{w^2|\nabla d_k|}{d_k^2(y) (\ln (j))^2} d_k^{2\alpha_k}(y) dy\le
C \frac{||w||_{L^\infty(\Omega)}^2}{(\ln(j))^2}
\int_{\frac{1}{j^2}<t< \frac{1}{j}} t^{2\alpha_k-2} t^{k-1} dt\le
$$
$$\le C \frac{||w||_{L^\infty(\Omega)}^2}{(\ln(j))^2}
\frac{1}{j^{2\alpha_k-2+k}} \int_{\frac{1}{j^2}<t< \frac{1}{j}}
t^{-1} dt \le C \frac{||w||_{L^\infty(\Omega)}^2}{\ln(j)}\to 0 \ ,$$
as $j\to \infty$ for any $2\alpha_k-2+k\ge 0$, and this completes
the proof.

\finedim

\bigskip
At this point we have all the ingredients needed in order to apply
Moser iteration technique up to the boundary, as adapted on
bounded domains in \cite{FMT1}, to the operator

\be \label{enne}{\mathcal L}_\alpha v:=-\frac{1}{\prod^n_{k=1}
d_k^{2\alpha_k}(x) } \sum^n_{i,j=1} \frac{\partial}{\partial x_i}
\left(a_{ij} (x) \ \prod^n_{k=1} d_k^{2\alpha_k}(x) \frac{\partial
v}{\partial x_j} \right)+\lambda_1 v \ ,\ee or equivalently to the
degenerate elliptic operator ${\mathcal L}_{\phi_1}$, defined in
(\ref{deflp}).

In fact one can prove the following result

\begin{theorem}{\bf (Parabolic Harnack inequality up to the
boundary)} \label{harnackgen} Let $n\ge 2$, $\Omega\subset \R^n$, be
a smooth bounded domain, $\lambda_1\in \R$ and $\alpha_k \ge
-\frac{k-2}{2}$, for $k=1,\cdots, n$. Then there exist positive
constants $C_H$ and $R=R(\Omega)$ such that for $x\in \Omega$,
$0<r<R$ and for any positive solution $v(y,t)$ of \be
\label{reference}\frac{\partial v}{\partial t}=-{\mathcal L}_\alpha
v \hbox{ in } \left\{{\mathcal B}(x,r)\cap \Omega\right\}\times
(0,r^2)\ ,\ee the following estimate holds true
$${\rm ess~sup}_{(y,t)\in \left\{{\mathcal B}(x,\frac{r}{2})\cap \Omega
\right\}\times (\frac{r^2}{4},\frac{r^2}{2})} v(y,t) \le C_H ~
{\rm ess~inf}_{(y,t)\in \left\{{\mathcal B}(x,\frac{r}{2})\cap
\Omega \right\}\times (\frac{3}{4} r^2,r^2)} v(y,t)\ .$$
\end{theorem}

Here we use the following definition of solutions:

\begin{definition}
\label{sol} By a solution $v(y,t)$ to (\ref{reference}), we mean a
function $$v\in C^1\left((0,r^2); L^2\left({\mathcal B}(x,r)\cap
\Omega, \prod^n_{k=1} d_k^{2\alpha_k} (y) dy \right)\right) \cap
C^0\left((0,r^2); H^1\left({\mathcal B}(x,r)\cap \Omega,
\prod^n_{k=1} d_k^{2\alpha_k} (y) dy\right)\right)$$ such that for
any $\Phi\in C^0((0,r^2);C^\infty_0({\mathcal B}(x,r) \cap \Omega))$
and any $0<t_1<t_2<r^2$ we have \be \label{ultima}
\int_{t_1}^{t_2}\int_{{\mathcal B}(x,r)\cap \Omega} \left\{ v_t
\Phi+ \sum_{i,j=1}^{n}  a_{ij}(y) \frac{\partial v}{\partial
x_i}\frac{\partial \Phi}{\partial x_j}+\lambda_1 v\ \Phi\right\}
\prod^n_{k=1} d_k^{2\alpha_k}(y) dy dt = 0 \ .\ee
\end{definition}

Let us note that Theorem \ref{harnackgen} is sharp, in the sense
that the same statement does not hold true if $\alpha_k <
-\frac{k-2}{2}$ for some $k=1,\ldots,n$ as explained in
\cite{FMT1}.

\smallskip
The parabolic Harnack inequality up to the boundary for the
Schr\"odinger type operator $\mathcal L$ defined in (\ref{defl})
and stated in Theorem \ref{harnackgen2} is proved as follows:

\smallskip
\noindent {\em Proof of Theorem \ref{harnackgen2}:} Clearly
Theorem \ref{harnackgen} applies also to the operator ${\mathcal
L}_{\phi_1}$ instead of ${\mathcal L}_{\alpha}$. Hence Theorem
\ref{harnackgen2} is a consequence of Theorem \ref{harnackgen} for
${\mathcal L}_{\phi_1}$ and of the change of variables
$v=u\phi_1^{-1}$ see ({\ref{change}).

\finedim

>From the parabolic Harnack inequality in Theorem \ref{harnackgen} we
deduce, as in \cite{FMT1}, the sharp two-sided estimates for the
heat kernel $l_\alpha$ associated to the elliptic operator
${\mathcal L}_\alpha$ defined in (\ref{enne}) under Dirichlet
boundary conditions. That is
$$v(x,t):=\int_{\Omega} l_\alpha(t,x,y) v_0(y) \prod^n_{i=1}
d_i^{2\alpha_i}(y)  dy$$ satisfies $v_t=-{\mathcal L}_\alpha v$ in
$ (0,\infty)\times \Omega $, $v(0,x)=v_0(x)$ on $\Omega$ and $v=0$
on $(0, \infty) \times \partial_1 \Omega $. We then have

\begin{theorem}
\label{heatgen} Let $n\ge 2$, $\Omega\subset \R^n$, be a smooth
bounded domain, $\lambda_1\in \R$ and $\alpha_k \ge -\frac{k-2}{2}$,
for $k=1,\cdots, n$. Then there exist positive constants $C_1, C_2$,
with $C_1\le C_2$, and $T>0$ depending on $\Omega$ such that
$$C_1 \prod^{n}_{i=1} (d_i(x)+\sqrt t)^{-\alpha_i}(d_i(y)+\sqrt t)^{-\alpha_i}
t^{-\frac{n}{2}} e^{-C_2\frac{|x-y|^2}{t}}\le l_\alpha(t,x,y) \le
$$ $$\le C_2 \prod^{n}_{i=1} (d_i(x)+\sqrt t)^{-\alpha_i}(d_i(y)+\sqrt t)^{-\alpha_i}
t^{-\frac{n}{2}} e^{-C_1\frac{|x-y|^2}{t}}$$ for all $x,y\in \Omega$
and $0<t\le T$.
\end{theorem}

Finally from the global upper bound in (\ref{uppp}), arguing as in
Theorem 6 of \cite{D1} (see also Proposition 4 in \cite{D2} as well
as the proof of Theorem 1.2 in \cite{FMT1}), one can deduce an
analogous lower bound for large times, thus obtaining the following
sharp long-time asymptotics of the heat kernel

\begin{theorem}
\label{heatgenlarge} Let  $n\ge 2$, $\Omega\subset \R^n$, be a
smooth bounded domain, $\lambda_1\in \R$ and $\alpha_k \ge
-\frac{k-2}{2}$, for $k=1,\cdots, n$. Then there exist positive
constants $C_1, C_2$, with $C_1\le C_2$, and $T>0$ depending on
$\Omega$ such that
$$C_1
e^{-\lambda_1 t} \le l_\alpha(t,x,y) \le C_2 e^{-\lambda_1 t}$$
for all $x,y\in \Omega$ and $t\ge T$.
\end{theorem}

>From Theorems \ref{heatgen} and \ref{heatgenlarge}, making use of
the equivalence (\ref{changebis}) as well as of assumption
(\ref{2.3x}), we get the corresponding result for the Schr\"odinger
operator $\mathcal L$ stated in Theorems \ref{heatgen2} in the
Introduction. We omit further details.

As we have already mentioned  integrating the sharp two-sided estimates for
$h(t,x,y)$ in Theorem \ref{heatgen2} with respect to the time
variable, one can deduced estimates on the Green function for the
Schr\"odinger operator ${\mathcal L}$ defined in (\ref{defl}) in the
case $\lambda_1>0$. Some explicit examples of sharp two sided Green function estimates
are given in Theorem 4.11 in \cite{FMT1}.

\section{Applications}
\setcounter{equation}{0}

In this Section we give some examples of singular potentials $V$ for
which the results of the present work apply; that is, we give
examples of potentials $V$ for which the generalized first
eigenvalue is not $-\infty$, and the corresponding first
eigenfunction is bounded from above and below uniformly by some
power of the distance function. We should stress that the
asymptotics of $\phi_1$ for the examples that follow  is a
consequence only of the maximum principle as used in \cite{BMS}. We
will present the detailed argument  for example III in the Appendix;
the other cases can be treated similarly.

To this end let us first prove that the sum of two potentials
having disjoint singularity sets and finite generalized first
eigenvalues, also has finite generalized first eigenvalue.

\medskip

\begin{lemma}\label{somma}
Let $V_i$, $i=1,2$ be such that $V_i\in L^1_{loc}(\Omega)\cap
L^\infty_{loc}(\overline{\Omega} \setminus S_i)$, where $S_i$ are
compact subsets of $\overline{\Omega}$ such that $S_1\cap
S_2=\emptyset$ and $$\lambda_1(V_i):=\inf_{u\in
C^\infty_0(\Omega)} \frac{\int_{\Omega} \left(|\nabla u|^2 -V_i
u^2\right) dx}{\int_{\Omega} u^2 dx} \ .$$ Assuming that
$\lambda_1(V_i)>-\infty$ for $i=1,2$ then
$\lambda_1(V_1+V_2)>-\infty$.
\end{lemma}

\smallskip
\noindent {\em Proof: } Let  us take $\varphi\in
C^\infty(\overline{\Omega})$, $0\le \varphi\le 1$, such that
$\varphi=1$ in $ \overline{\Omega}\cap\Omega_1$ and $\varphi=0$ in
$\overline{\Omega} \setminus \tilde\Omega_1$ where $\Omega_1$ is a
neighborhood of $S_1$, that is $\Omega_1:=\{x\in
\overline{\Omega}: {\rm dist} (x,S_1)<\delta \}$ for some small
$\delta>0$, and $\tilde \Omega_1$ is a slightly bigger
neighborhood of $S_1$, thus $ S_1 \subset \Omega_1\subset \tilde \Omega_1$.
Whence $\Omega_2:=\overline{\Omega}\setminus \tilde\Omega_1$ is a
neighborhood of $S_2$ and $\tilde \Omega_2:=\overline{\Omega}
\setminus \Omega_1$ is a slightly bigger neighborhood of $S_2$,
thus $\Omega_2\subset \tilde \Omega_2$. Whence for any $u\in
C^\infty_0(\Omega)$ we have $u=u\varphi+u(1-\varphi)=:u_1+u_2$. By
elementary calculations we have that
$$\int_{\Omega}|\nabla u|^2dx=\int_{\Omega} |\nabla
(u\varphi+(1-\varphi) u)|^2 dx \ge \int_{\Omega}|\nabla u_1|^2dx
+\int_{\Omega}|\nabla u_2|^2 dx -K\int_{\Omega}u^2 dx$$ for a
suitable positive constant $K$. Then for $V:=V_1+V_2$ we have
$$\int_{\Omega} \left(|\nabla u|^2-V u^2 \right) dx \ge $$
$$\ge\int_{\Omega} \left(|\nabla u_1|^2-V_1 u_1^2 + |\nabla
u_2|^2-V_2u_2^2 + u^2(V_1\varphi^2+V_2(1-\varphi)^2-V-K)\right)dx
=$$ $$= \int_{\tilde \Omega_1}\left( |\nabla u_1|^2-V_1 u_1^2
\right) dx+\int_{\tilde \Omega_2}\left(|\nabla u_2|^2-V_2
u_2^2\right) dx+$$
$$+\int_{\Omega}u^2\left[V_1(\varphi^2-1)+V_2((1-\varphi)^2-1)-K\right]
dx \ge$$
$$\ge
\left(\lambda_{1}(V_1)+\lambda_{1}(V_2)-||V_1||_{L^\infty(\tilde
\Omega_2)} -||V_2||_{L^\infty(\tilde
\Omega_1)}-K\right)\int_{\Omega} u^2 dx \ .$$

\finedim

We now present some concrete examples.

\smallskip
\noindent {\bf Example I } Our first example is motivated by
\cite{FeT1}, \cite{FeT2}, \cite{FeMaT} and deals with multipolar
inverse-square potentials. Let $n\ge 3$, $\Omega\subset \R^n$ be a
 smooth bounded domain from which we have removed $m$ points $x_1,\ldots,x_m$ and
$$V(x)=\sum^m_{i=1}\frac{c_i}{|x-x_i|^2} \ ,$$

\noindent for $0\le c_i \le \frac{(n-2)^2}{4}$. We note that
differently from \cite{FeT1}, we may take in each one of the
inverse-square potentials the critical Hardy constant. This is due
to the fact that we study the Schr\"odinger operator $-\Delta -V$ on
a bounded domain. In such a case one can prove that $\phi_1(x)\sim
\prod^{m}_{i=1} |x-x_i|^{\beta_i} {\rm dist }(x,\partial_1 \Omega)$
with
$$\beta_i:=\frac{2-n+\sqrt{(n-2)^2-4c_i}}{2}  \ ,$$ see Lemma 7 in
\cite{DD2} and Theorem 7.1 in \cite{DS} on one hand and the
elliptic regularity on the other. In fact the function
$f(x):=|x-x_i|^{\beta_i}$ satisfies the equation $\Delta
f+\frac{c_i}{|x-x_i|^2}f=0$ in $\R^n\setminus\{x_i\}$. We only
need to check that $\lambda_1>-\infty$. This follows from Lemma
\ref{somma}, which clearly can be generalized to a finite sum of
potentials, and the improved $L^2$ inequality given in \cite{VZ}
for a single inverse-square potential. Consequently Theorems
\ref{harnackgen2}, \ref{heatgen2} and \ref{thm3.1} apply with
$\Gamma_k=\emptyset$ for $k=2,\ldots,n-1$ and
$\Gamma_n=\{x_1,\ldots,x_m\}$. Note that in this case
$d_n^{\alpha_n}(x)$ stands for $\prod^m_{i=1} |x-x_i|^{\beta_i}$,
that is, we have $m$ different sets of the same codimension $n$ in
the boundary of $\Omega$, where $\phi_1$ may present different
degeneracies.

\smallskip
\noindent {\bf Example II } Let $n\ge 4$, $\Omega=B_R\setminus E$,
for some $R>1$, where $$E:=\{x\in \R^n: x_1^2+x_2^2=1, \
x_3=\ldots=x_n= 0\}$$
$$B_R:=\{x\in \R^n: x_1^2+\ldots +x_n^2 < R^2\} \ \hbox { and }$$
$$V(x)=\frac{1}{4}\frac{1}{{\rm dist}^2(x,\partial
B_R)}+\frac{(n-3)^2}{4}\frac{1}{{\rm dist}^2(x,E)} \ .$$  In such a
case one can easily prove that $\phi_1(x)\sim {\rm
dist}^{\frac{1}{2}}(x,\partial B_R)\ {\rm dist}^{\frac{3-n}{2}}(x,
E)$, see \cite{DD2}. The fact that $\lambda_1>-\infty$ follows
making use of Lemma \ref{somma}, from
 the improved $L^2$ inequality given in
\cite{BM} for the inverse-square distance to $\partial B_R$  and the
one given in \cite{DD2} or \cite{FMaT2} for the inverse-square
distance to the set $E$ having codimension $n-1$. Theorems
\ref{harnackgen2}, \ref{heatgen2} and \ref{thm3.1} now apply to
$-\Delta -V$ with $\alpha_1=1/2$ and $\alpha_{n-1}=(3-n)/2$, whereas
all the other $\alpha_k$'s are zero.

\smallskip
\noindent {\bf Example III } Let $n\ge 2$, and  $\Omega\subset \R^n$
be a bounded domain such that $\partial \Omega = \partial_1 \Omega$,
that is, the boundary of $\Omega$ has codimension one. We now take
$$V(x)=\frac{1}{4}\frac{1}{{\rm dist}^{2}(x,\partial \Omega)} \
.$$ \noindent
By the results of \cite{BM}
 we  have that
$\lambda_1>-\infty$ under appropriate regularity assumptions on
$\partial \Omega$.  We recall also that $\phi_1(x)\sim{\rm
dist}^{\frac{1}{2}}(x,\partial \Omega)$, as shown in \cite{DD2}; see
also Appendix, where we will provide a self-contained proof based
only on the   maximum principle. Therefore Theorems
\ref{harnackgen2}, \ref{heatgen2} and \ref{thm3.1} apply to $-\Delta
-V$ with $\alpha_1=1/2$, whereas all the other $\alpha_k$'s are
zero. We note that this improves the corresponding Theorems in
[FMT1] removing the convexity assumption, under which it is known
that $\lambda_1>0$ (see \cite{BM}).

\smallskip
\noindent {\bf Example IV } Let $n\ge 3$, $\Omega' \subset \R^n$ be
a smooth  bounded domain containing the origin,  $\Omega = \Omega' \setminus \{0\}$,  and
$$V(x)=\frac{1}{4}\frac{1}{{\rm dist}^2(x,\partial_1
\Omega)}+\frac{(n-2)^2}{4}\frac{1}{|x|^2}\ .$$ \noindent The fact
that $\lambda_1>-\infty$ may be deduced from Lemma \ref{somma}
making use of the $L^2$ improved Hardy inequality in \cite{VZ} for
the inverse-square potential $1/|x|^2$  and of the one  in \cite{BM}
for the inverse-square potential involving the distance to
$\partial_1 \Omega$. In this example we have $\phi_1(x)\sim {\rm
dist}^{\frac{1}{2}}(x,\partial_1 \Omega)|x|^{\frac{2-n}{2}}$. Whence
Theorems \ref{harnackgen2}, \ref{heatgen2} and \ref{thm3.1} apply to
$-\Delta-V$ with $\alpha_1=1/2$ and $\alpha_{n}=(2-n)/2$, whereas
all the other $\alpha_k$'s are zero.

\smallskip
\noindent {\bf Example V }  Let $n\ge 3$, and $B_1 \subset \R^n$ be the unit ball.
For $-\frac{n-2}{2} \leq a < 0$ we
 consider the operator  ${\mathcal L} = -L -V$ where
\[
Lu = \sum^{n}_{i,j=1} \frac{\partial}{\partial x_i}
\left(a_{ij}(x) \frac{\partial
u}{\partial x_j}\right),~~~~~~~a_{ij}=\delta_{ij} +\frac12 |x|^{2-a}(1-\delta_{ij}),
\]
and
  $V(x) = -\frac{a(n+a-2)}{|x|^2}$.  The operator $L$ is easily seen to be uniformly elliptic
and in fact
\be\la{elip}
\left( 1 + \frac{n}{2} \right) |\xi|^2 \geq  \sum^{n}_{i,j=1}a_{ij} \xi_i \xi_j \geq \frac12 |\xi|^2.
\ee
On the other hand if
\[
Q[u]= \int_{B_1} \left( \sum_{i,j=1}^{n} a_{ij}(x) \frac{\partial
u}{\partial x_i}\frac{\partial u}{\partial x_j} - V u^2 \right)
dx ,
\]
using the change of variables $u= |x|^a v$ a straightforward calculation shows that
\[
 \int_{B_1} \left(\sum_{i,j=1}^{n} a_{ij}(x) \frac{\partial
u}{\partial x_i}\frac{\partial u}{\partial x_j} - V u^2\right)
dx =  \int_{B_1} (\delta_{ij} |x|^{2a} + \frac12 |x|^{a+2} (1-\delta_{ij})) v_{x_i}v_{x_j} dx .
\]
Using (\ref{elip}) one can easily  see that $\lambda_1>0$.

In this example we have that $\phi_1(x) \sim {\rm dist}(x,\partial
B_1) |x|^{a}$. We can apply our results with  $\alpha_1=1$ and
$\alpha_{n}=a$, whereas all the other $\alpha_k$'s are zero.

\section{Appendix}
\setcounter{equation}{0}

In this Appendix we consider the operator ${\mathcal L}:=-\Delta
-\frac{1}{4}\frac{1}{d^2(x)}$ which corresponds to Example III of
Section 4 and we will prove that the corresponding eigenfunction
$\phi_1$ is such that $\phi_1(x) \sim d^{\frac{1}{2}}(x)$, where
$d(x):={\rm dist}(x,\partial \Omega)$.

\smallskip
\noindent  {\it Step I: Existence of $\phi_1$ in a suitable energy space}.
Let $\eta \in C^2(\Omega)$    be a function such that $\eta(x)=d^{1/2} (x)$ near the boundary,
say,   $d(x) \leq \epsilon_0$, and $ \eta(x) \geq c_0 >0$  for $d(x) \geq \epsilon_0$.
 Let $H^1_d(\Omega)$ be the
closure of $C^\infty_0(\Omega)$ functions
 under the norm
$$||v||^2_{H^1_{d}}:=\int_{\Omega} d \left(|\nabla v|^2 +v^2
\right) dx  \ .$$
This norm is equivalent to the norm
$
\|v\|:=\int_{\Omega} \eta^2 \left(|\nabla v|^2 +v^2
\right) dx.
$
Changing  variables by  $u = \eta v$, in
\[
- \infty  < \lambda_1 = \inf_{u \in C^\infty_0(\Omega)} \frac{\int_{\Omega}( |\nabla u|^2 - \frac{u^2}{4 d^2}) dx}{\int_{\Omega} u^2 dx},
\]
we get  the equivalent inequality
\be\label{ap1}
- \infty  < \lambda_1  = \inf_{v \in C^\infty_0(\Omega)}
 \frac{\int_{\Omega} (\eta^2 |\nabla v|^2 - (\eta \Delta \eta + \frac{\eta^2}{4 d}) v^2)dx}
{\int_{\Omega} \eta^2 v^2 dx}.
\ee

Using   the fact  that  $  \eta \Delta \eta + \frac{\eta^2}{4 d}  \in L^{\infty}(\Omega)$
 as well  as the following estimate
\[
\int_{\Omega} \frac{X^2(d)}{d}
v^2 d x \le C \int_{\Omega} d (  |\nabla v|^2 +  v^2)  dx    \ ,
 \ \
~~~ \ v\in C^\infty_0(\Omega),
\]
which was established in Proposition 5.1 in \cite{FMT2}, we find
that for every $\epsilon >0$ there exists an $M_{\epsilon}>0$ such
that for all $v\in C^\infty_0(\Omega)$, \be\label{ap2} \Big|
\int_{\Omega} (\eta \Delta \eta + \frac{\eta^2}{4 d})
 v^2 dx \Big| \leq \epsilon \int_{\Omega} \eta^2  |\nabla v|^2 dx  + M_{\epsilon} \int_{\Omega}  \eta^2 v^2 dx.
\ee
In the sequel we will establish the existence of a function $\psi_1 \in H^1_{d}(\Omega)$ which realizes the infimum in
(\ref{ap1}).  To this  end let $w_k$ be a minimizing sequence normalized by  $ \int_{\Omega}  \eta^2 w_k^2 dx=1$.
Then, using (\ref{ap2})  we can easily obtain that  the sequence $w_k$ is bounded in $H_d^1$. Therefore there
exists a subsequence still denoted by $w_k$ such that it converges  $H^1_d$ -- weakly to $\psi_1$, and in addition we
have
the following strong convergence, for $k \rightarrow \infty$,
\[
 \int_{\Omega} (\eta \Delta \eta + \frac{\eta^2}{4 d})
 w_k^2 dx  \rightarrow  \int_{\Omega} (\eta \Delta \eta + \frac{\eta^2}{4 d})
 \psi_1^2 dx,
\]
and
\[
 \int_{\Omega} \eta^2 w_k^2 dx  \rightarrow  \int_{\Omega}  \eta^2  \psi_1^2 dx.
\]
Using the lower semicontinuity of the gradient term in the numerator of (\ref{ap1}) the result
follows.

\noindent {\it Step II: An auxiliary estimate.} For $\delta>0$ small
enough we  set \be \label{lastde} \mu_1(\Omega_{\delta}):=
\inf_{v\in C^\infty_0(\Omega_\delta)} \frac{\int_{\Omega_\delta}
\left(d |\nabla v|^2 -\Delta d \ \frac{v^2}{2}\right)
dx}{\int_{\Omega_\delta}d  v^2 dx} \ ,\ee where
$\Omega_\delta=\{x\in \Omega \hbox{ s.t. dist }(x,\partial
\Omega)<\delta\}.$

We will show that \be\label{ap3} \mu_1(\Omega_{\delta})  \to +
\infty, ~~~~~~~~~~{\rm as} ~~~~~~~~ \delta \to 0. \ee Our starting
point is the inequality \be \label{disn}\int_{\Omega_\delta}
\left(|\nabla u|^2-\frac{1}{4}\frac{u^2}{d^2}\right) dx \ge
\frac{1}{8} \int_{\Omega_\delta} \frac{X^2(d)}{d^2} u^2 dx \ , \ \
~~~ \ u \in C^\infty_0(\Omega_\delta),\ee for any $0<\delta \le
\delta_0$, for some $\delta_0$ small enough, where $X(t):=(1-\ln
t)^{-1}$. To prove this  one starts with  the obvious relation
$$0\le \int_{\Omega_{\delta}} \Big |\nabla u -\left(\frac{\nabla
d}{2d}-\frac{X\nabla d}{2d}\right) u\Big|^2 dx \ .$$ Expanding the
square, integrating by parts and using the fact that $|d\Delta d|$
can be made arbitrarily small in $\Omega_\delta$, for $\delta$
sufficiently small, the result follows. Changing variable as usual
by $u=d^{\frac{1}{2}}v$, inequality (\ref{disn}) is equivalent to
\be \label{reallast}\int_{\Omega_\delta} \left(d|\nabla
v|^2-\Delta d \frac{v^2}{2}\right) dx \ge \frac{1}{8}
\int_{\Omega_\delta} \frac{X^2(d)}{d} v^2 dx \ , \ \ ~~~ \ v
\in C^\infty_0(\Omega_\delta)\ .\ee
Therefore
$$\frac{1}{8} \frac{X^2(\delta)}{\delta^2}\le
  \frac{X^2(\delta)}{\delta^2}\frac{\int_{\Omega_\delta} \left(d|\nabla v|^2-\Delta d
\frac{v^2}{2}\right) dx}{\int_{\Omega_\delta} \frac{X^2(d)}{d} v^2
dx}\le\frac{\int_{\Omega_\delta} \left(d|\nabla v|^2-\Delta d
\frac{v^2}{2}\right) dx}{\int_{\Omega_\delta}  v^2 d \ dx}, $$
from which   (\ref{ap3}) follows.

\noindent
{\it Step III: Asymptotics of $\phi_1$.} The lower bound $C_1 d^{1/2}(x) \leq \phi_1(x)$ is a consequence
of the maximum principle and is derived in
 Lemma 7 in
\cite{DD2}.

 We will  obtain the upper bound using  maximum principle  in a suitably small neighborhood
of the boundary.  Let $\psi_1(x)= \phi_1(x)/d^{1/2}$. Then, for $E
v:=-{\rm div}(d \nabla v) -\frac{\Delta d}{2}v - \lambda_1 d v$, we
have that $E\psi_1 = 0$. Moreover, we have that
\[
E(1-Cd)=C-\frac{\Delta d}{2} +d \left(-\lambda_1+\frac{3C}{2} \Delta d
+C \lambda_1  d \right) \geq 0,
\]
in $\Omega_{\delta}$ for $\delta$ small enough and $C>0$ big enough. We next choose $\beta>0$ big enough
so that
\[
\psi_1(x)  \leq \beta (1-Cd)~~~~~~~~~{\rm on} ~~~~~~~\partial \Omega_{\delta}.
\]
Let $g(x) := \psi_1(x) - \beta(1-Cd)$  and $g^{+}:=\max\{0,g\}$.  We clearly have that
\[
 \int_{\Omega_\delta} g^+ Eg \leq 0,
\]
from which it follows that
\[
\frac{\int_{\Omega_\delta}( d |\nabla g^+|^2 -\frac{\Delta d}{2}(g^{+})^2 ) dx}{\int_{\Omega_\delta} d (g^{+})^2 dx} \leq
\lambda_1.
\]
This contradicts   (\ref{ap3})  unless $g^+=0$ from which it
follows that $\phi(x) \leq \beta d^{1/2}(x)$.

\medskip
\noindent {\bf Acknowledgments} LM acknowledges the support of
University of Crete and FORTH as well as the "Progetto di Ateneo
Federato" 2007, Universit\'a di Roma La Sapienza, "Elliptic and
parabolic equations, minimum of functionals: existence of solutions
and qualitative properties", during her visits to Greece. AT
acknowledges the support of Universities of Rome I, Bologna and
FORTH as well as the GNAMPA project "Liouville theorems in
Riemannian and sub-Riemannian settings" during his visits in Italy.

The authors thank the referee  for his comments and suggestions.

\end{document}